\documentclass[10pt,a4paper,reqno]{amsart}
\usepackage{amssymb}
\usepackage{amscd}
\usepackage{hyperref}
\numberwithin{equation}{section}

\theoremstyle{plain}

\newtheorem{theorem}[subsection]{Theorem}
\newtheorem{proposition}[subsection]{Proposition}
\newtheorem{lemma}[subsection]{Lemma}
\newtheorem{corollary}[subsection]{Corollary}

\theoremstyle{definition}

\newtheorem{remark}[subsection]{Remark}
\newtheorem{remarks}[subsection]{Remarks}
\newtheorem{example}[subsection]{Example}
\newtheorem{examples}[subsection]{Examples}

\renewcommand{\leq}{\leqslant}
\renewcommand{\geq}{\geqslant}

\newsavebox{\proofbox}
\savebox{\proofbox}{\begin{picture}(7,7)%
  \put(0,0){\framebox(7,7){}}\end{picture}}

\newcommand\Z{\mathbb{Z}}
\newcommand\R{\mathbb{R}}
\newcommand\T{\mathbb{T}}

\newcommand\N{\mathbb{N}}

\def\proof{\noindent\textit{Proof. }}

\def\remark{\noindent\textit{Remark. }}

\def\endproof{\hfill{\usebox{\proofbox}}}

\begin{document}

\title[Bilinear Estimates]{A bilinear oscillatory integral estimate and bilinear refinements to Strichartz estimates on closed manifolds}
\author{Zaher Hani}

\address{UCLA Department of Mathematics, Los Angeles, CA 90095-1555.}
\email{zhani@math.ucla.edu}

\maketitle

\begin{abstract}We prove a bilinear $L^2(\R^d) \times L^2(\R^d) \to L^2(\R^{d+1})$ estimate for a pair of oscillatory integral operators with different asymptotic parameters and phase functions satisfying a transversality condition. This is then used to prove a bilinear refinement to Strichartz estimates on closed manifolds, similar to that derived in \cite{B2} on $\R^d$, but at a relevant semi-classical scale. These estimates will be employed elsewhere \cite{H} to prove global well-posedness below $H^1$ for the cubic nonlinear Schr\"odinger equation on closed surfaces.
\end{abstract}

\section{Introduction}

We consider oscillatory integrals defined by:

\begin{equation}\label{T_lambda}
T_\lambda f(t,x)=\int_{\R^d} e^{i\lambda\phi(t,x,\xi)}a(t,x,\xi)f(\xi)d\xi
\end{equation}

where $t\in \R$, $x,\xi \in \R^d$, $a\in C_0^\infty(\R\times \R^d \times \R^d)$. The phase function $\phi$ is a real-valued smooth function on the support of $a$. We shall assume that it satisfies a usual non-degeneracy condition, namely that the $(n+1)\times n$ matrix:

\begin{equation}\label{nondegeneracy}
\frac{\partial^2 \phi}{\partial \xi \partial(x,t)}(t_0,x_0,\xi_0) \textrm{ has maximal rank $d$ for every $(t_0,x_0,\xi_0)\in \operatorname{supp} a$}
\end{equation}

This implies that for each fixed $(t_0,x_0)\in \R^{d+1}$, the map given by:

$$
\xi \mapsto \nabla_{(t,x)}\phi(t_0,x_0,\xi)
$$ 

defines a smooth immersion from $\R^d$ into $\R^{d+1}$. The image of this map is a hyper-surface which we denote by $S_\phi(t_0,x_0)$ and $S_\phi$ when no confusion arises. Our objective is to prove bilinear estimates for such operators and use them to get bilinear refinements to Strichartz estimates on compact manifolds without boundary.

Operators as in $\eqref{T_lambda}$ can be thought of as variable coefficient generalizations of usual dual restriction (extension) operators where $\phi(t,x,\xi)=x.\xi +t\psi(\xi)$ and $\eqref{T_lambda}$ becomes the dual of the operator given by restricting the Fourier transform to the hyper-surface $S_\phi=\{(\tau,\xi)\in \R^{d+1}:\tau=\psi(\xi)$\}. As in the case of restriction operators, one is interested in obtaining asymptotic decay estimates for $||T_\lambda||_{L^p(\R^d)\to L^q(\R^{d+1})}$ in terms of $\lambda$. It is well known that in order to obtain non-trivial decay estimates (the optimal one being $\lambda^{-\frac{d+1}{q}}$), one has to impose some curvature condition on the hyper-surfaces $S_\phi$, namely that the Gaussian curvature does not vanish anywhere. The pairs of exponents $(p,q)$ for which this decay is possible were specified by H\"ormander in \cite{H1} when $d=1$ and posed as a question for higher dimensions. Since then, there has been a tremendous amount of research in proving such bounds. (see \cite{Stein} and references therein for an introduction and \cite{TPC} for a more current survey). 

We will be interested in bilinear versions of such estimates. In this case, one considers the product $T_\lambda f \tilde T_\mu g$ where $\tilde T_\mu g$ is an operator similar to $\eqref{T_lambda}$

\begin{equation}\label{tilde T_mu}
\tilde T_\mu g(t,x)=\int_{\R^d} e^{i\mu\psi(t,x,\xi)}b(t,x,\xi)g(\xi)d\xi
\end{equation}

where $b\in C_0^\infty(\R\times \R^d \times \R^d)$ and $\psi$ is smooth on the support of $b$ and satisfies the same non-degeneracy assumption $\eqref{nondegeneracy}$. The initial motivation behind such estimates was proving and refining the linear estimates in the case when the exponent $q$ is an even number. However, such an improvement is only possible when the surfaces $S_\phi$ and $S_\psi$ satisfy a certain transversality assumption. This transversality turns out to be more important than any curvature assumption in certain instances. To be precise, the type of estimates one is often interested in are of the form:

\begin{equation}\label{bilinear restriction}
||T_\lambda f \tilde T_\mu g||_{L^q(\R\times \R^d)}\lesssim \Lambda(\lambda,\mu)||f||_{L^2(\R^d)}||g||_{L^2(\R^d)}
\end{equation}

(For us, the case when $q=2$ and $\lambda  \neq \mu$ will be of particular interest). Great progress has been achieved in proving estimates like $\eqref{bilinear restriction}$ especially in the case $\lambda=\mu$ and when the surfaces $S_\phi$ and $S_\psi$ satisfy some non-vanishing curvature assumption. In the constant coefficient (restriction) case, Wolff was able to prove $\eqref{bilinear restriction}$ in the cone restriction case for all $q>1+\frac{2}{d+1}$ with $\Lambda(\lambda, \lambda)\lesssim \lambda^{-\frac{d+1}{q}}$\cite{W} . This estimate was later extended to the endpoint by Tao in \cite{Tao endpoint}. The same estimate was then proven for transverse subsets of the paraboloid \cite{Tao paraboloids}. In the variable coefficient case, Lee proved a similar estimate when $\lambda=\mu$, $q\geq 1+\frac{2}{d+1}$, and $\Lambda(\lambda, \lambda)\lesssim \lambda^{-\frac{d+1}{q}+\epsilon}$ under certain curvature assumptions on the surfaces $S_\phi(t_0,x_0)$ and $S_\psi(t_0,x_0)$ \cite{Lee}. 

In this paper, we prove an $L^2$ estimate when $\lambda \neq \mu$ and the only assumption we impose on the hyper-surfaces $S_\phi$ and $S_\psi$ is transversality. In particular, no curvature assumptions are taken. 

\begin{theorem}\label{Bilinear FIO estimate}
Suppose that $T_\lambda$ and $\tilde T_\mu$ are two oscillatory integral operators of the form given in $\eqref{T_lambda}$ with $\mu \leq \lambda$ and assume that the canonical hyper-surfaces associated with the phase functions $\phi$ and $\psi$ satisfy the standard transversality condition $\eqref{transversality condition}$, then:

\begin{equation}\label{bilinear estimate}
\left|\left|T_{\lambda}f \tilde T_\mu g\right|\right|_{L^2(\R\times \R^d)}\lesssim \frac{1}{\lambda^{d/2}\mu^{1/2}}||f||_{L^2(\R^d)}||g||_{L^2(\R^d)}.
\end{equation}

The implicit constants are allowed to depend on $\delta$, $d$, and uniform bounds on a fixed number of derivatives of $\phi, \psi, a$, and $b$.

\end{theorem}

A couple of remarks are in order. First, we mention that $\eqref{bilinear estimate}$ is sharp (cf. remark at the end of section \ref{proof of main theorem}). Second, we note that without curvature assumptions on the surfaces, the linear estimate is easily seen to fail (e.g. restriction to hyperplanes). However, the $L^2$ bilinear estimate is true as long as the surfaces are transverse.\footnote{This is well-known in the constant coefficient case, see \cite{TPC}.}. Even when the linear estimate is true (which requires as mentioned a non-vanishing curvature assumption on the surfaces), $\eqref{bilinear estimate}$ is an improvement on applying H\"older and the linear estimates available especially in the case when $\mu \ll \lambda$ (for example, when $d=2$ linear estimates give the bound$\frac{1}{\lambda^{3/4}\mu^{3/4}}$). This improvement is often of great importance in applications (see \cite{B3},\cite{B2}, \cite{H}).

We now specify the transversality condition needed. The canonical hyper-surfaces $S_\phi(t_0,x_0)$ and $S_\psi(t_0,x_0)$, given by the maps $\xi \mapsto \nabla_{(t,x)}\phi(t_0,x_0,\xi)$ and $\xi \mapsto \nabla_{(t,x)}\psi(t_0,x_0,\xi)$ respectively, live in the cotangent space $T^*_{(t_0,x_0)}\R^n$ to $\R^n$ at $(t_0,x_0)$. The non-degeneracy condition defined in $\eqref{nondegeneracy}$ for $\phi$ (and defined similarly for $\psi$), implies that for every $\xi_0 \in \operatorname{supp}_\xi a(t_0,x_0,.)$, there exists a locally defined unit normal vector field $\nu_1(t_0,x_0,\xi_0)=\nu_1(\xi_0)$ to this surface at the point $\nabla_{(t,x)}\phi(t_0,x_0,\xi_0) \in T^*_{(t_0,x_0)}\R^{n+1}$. In other words, the map

$$
\xi \mapsto \langle \nu_1(\xi_0), \nabla_{(t,x)} \phi (t_0,x_0,\xi)\rangle 
$$

has a critical point at $\xi=\xi_0$ (in linear algebra terms, $\nu(\xi_0)$ is the unit vector spanning the one dimensional orthogonal complement of the image of the matrix appearing in $\eqref{nondegeneracy}$).
Similarly, we define the associated unit normal vector $\nu_2(\xi_0)$ to $S_\psi(t_0,x_0)$ at the point $\nabla_{(t,x)}\psi(t_0,x_0,\xi_0)$ satisfying:

$$
\xi \mapsto \langle \nu_2(\xi_0), \nabla_{(t,x)} \psi (t_0,x_0,\xi)\rangle 
$$

has a critical point at $\xi=\xi_0$.

The transversality condition we impose on the phase functions $\phi$ and $\psi$ is that \emph{the two surfaces $S_\phi(t_0,x_0)$ with $S_\psi(t_0,x_0)$ are uniformly transverse for every $(t_0,x_0)$}: by which we mean that there exists a $\delta>0$ such that for each $(t_0,x_0,\xi_1)\in \operatorname{supp } a$, $(t_0,x_0,\xi_2)\in \operatorname{supp }b$, we have:

\begin{equation}\label{transversality condition}
\left|\langle \nu_1(\xi_1),\nu_2(\xi_2)\rangle\right| \leq 1-\delta.
\end{equation}

This transversality condition is standard in all bilinear oscillatory integral estimates. We remark that there is a slight difference between this definition of transversality and that used in most differential topology textbooks in which the definition of transversality includes manifolds that do not intersect. Here we say that two hyper-surfaces are transverse if the intersection of all their translates is transverse in the sense of differential topology.

\remark The phase functions $\phi$ and $\psi$ can depend on $\lambda$ and $\mu$ as long as the quantitative estimates needed in the proof (namely $\eqref{transversality condition}$ and the derivative bounds mentioned in Theorem \ref{bilinear estimate}) are satisfied uniformly in $\lambda$ and $\mu$ on the support of $a$ and $b$.

The proof of Theorem \ref{Bilinear FIO estimate} is based on a $TT^*$ argument and delicate analysis of a cumulative phase function.

\subsection{Bilinear Strichartz Estimates} Our main application of the bilinear estimate in Theorem \ref{Bilinear FIO estimate} is to derive short-range or semi-classical bilinear Strichartz estimates for the Schrodinger equation on closed (compact without boundary) $d-$manifolds $M^d$. We will also be able to prove mixed bilinear estimates of Schr\"odinger-Wave type as well (see section \ref{further results}). Bilinear estimates are of great importance in PDE as they offer refinements to linear Strichartz estimates. The latter are given on $\R^d$ with its Euclidean Laplacian by:

\begin{equation}\label{linear SE on R^n}
||e^{it\Delta}u_0||_{L_t^qL_x^r(\R\times \R^d)}\lesssim ||u_0||_{L^2(\R^d)}
\end{equation}

where $(q,r)$ is any \emph{Schr\"odinger admissible} pair, i.e. $2\leq q,r \leq \infty$, $\frac{2}{q}+\frac{d}{r}=\frac{d}{2}$, and $(q,r,d)\neq (2,\infty,2)$. The implicit constants depend on $(q,r,d)$. These estimates are of fundamental importance in proving both local and global results for nonlinear Schr\"odinger equations. (cf. \cite{Tao book},\cite{KT}).

In the case of compact manifolds, the first Strichartz estimates were proved by Bourgain \cite{B1} in the case of the torus. The case of general compact Riemannian manifolds $(M,g)$ without boundary was dealt with by Burq, Gerard, and Tzvetkov in \cite{BGT} and \cite{ST}. In \cite{BGT}, the authors prove the following estimates:

\begin{equation}\label{BGT linear estimate}
||e^{it\Delta_g}u_0||_{L_t^q L_x^r([0,1]\times M)}\lesssim_{q,r,M} ||u_0||_{H^{\frac{1}{q}}(M)} 
\end{equation}

for any admissible pair $(q,r)$. The proof relies on a construction of an approximate parametrix to the semi-classical operator $e^{ih\Delta_g}\varphi(h\sqrt{-\Delta_g})$ (where $\varphi$ is Schwartz) which is used to prove the following \emph{semiclassical linear Strichartz estimate}:

\begin{equation}\label{semiclassical linear estimate}
||e^{it\Delta_g}u_0||_{L_t^qL_x^r([0,\frac{\alpha}{N}]\times M)}\lesssim_{q,r, M} ||u_0||_{L^2(M)}
\end{equation}

whenever $u_0$ is frequency (spectrally) localized at the dyadic scale $N$ and $\alpha \ll 1$. This estimate conforms with the heuristic that Schr\"odinger evolution moves wavepackets localized at frequency $\sim N$ at speeds $\sim N$, which means that in the time interval $[0,\frac{\alpha}{N}]$, one expects the wave packet to remain in a coordinate patch and hence satisfy the same estimates like those on $\R^d$. This heuristic will be very useful in predicting the right bilinear estimate later on as well.
Notice that $\eqref{BGT linear estimate}$ follows directly from $\eqref{semiclassical linear estimate}$ by splitting the time interval $[0,1]$ into $N$ subintervals of lengths $N^{-1}$ and using the conservation of mass and a square function estimate (cf. \cite{BGT}).

Turning to bilinear estimates, we will start by mentioning the relevant estimate on $\R^d$ for which we wish to find an analogue on compact manifolds. This estimate first appeared as a refinement to linear Strichartz estimates in Bourgain's paper \cite{B2}: assuming that $u_0$ is frequency localized at frequencies $\{\xi\in \R^d: |\xi|\sim N_1\}$ and $v_0$ is frequency localized at frequencies $\{\xi\in \R^d: |\xi|\lesssim N_2\}$ with $N_2 \leq N_1$, then the following holds:

\begin{equation}\label{BLS on R^d}
||e^{it\Delta}u_0 e^{it\Delta}v_0||_{L^2(\R \times \R^d)} \lesssim_d \frac{N_{2}^{\frac{d-1}{2}}}{N_{1}^{\frac{1}{2}}}||u||_{L^2(\R^d)}||v||_{L^2(\R^d)}
\end{equation}

We first notice that this estimate is an improvement on applying H\"older's inequality and the linear Strichartz estimates. In fact, applying the linear estimates only, one would get instead of the $\left(\frac{N_2^{(d-1)/2}}{N_1^{1/2}}\right)$ constant on the LHS of $\eqref{BLS on R^d}$: 1 for $d=2$ (here one uses the $L_x^2 \to L_{t,x}^4$ Strichartz estimate) and $N_2^{d/2-1}$ for $d\geq 3$ (here one should use H\"older, the $L_x^2 \to L_{t,x}^{\frac{2(d+2)}{d}}$ estimate for $e^{it\Delta}u_0$, and Bernstein combined with the $L_x^2 \to L_t^{d+2}L_x^{\frac{2d(d+2)}{d(d+2)-4}}$ for $e^{it\Delta} v_0$). Bourgain used this improvement (when $N_2 \ll N_1$) to prove, among other things, global well-posedness below energy norm for certain mass (and $\dot H^{1/2}$)-critical equations (which incidentally is also an application that will be considered in the context of closed manifolds in \cite{H}). Since then, this improvement and variants of it proved to be of essential use in studying nonlinear Schr\"dinger equations.

In the context of compact manifolds, some bilinear estimates on the torus were already implicit in the work of Bourgain \cite{B1}(cf. \cite{BEE}) and other variants were proved in \cite{dSPST}. In \cite{BEE} and \cite{MEE}, the authors prove bilinear Strichartz estimates on spheres $S^2$ and $S^3$ (and on the bit wider class of Zoll manifolds) using bilinear eigenfunction cluster estimates. These bilinear Strichartz estimates take the form:

$$
||e^{it\Delta_g}u_0 e^{it\Delta_g}v_0||_{L_{t,x}^2([0,1]\times S^d)}\lesssim_d N_2^{\alpha_d}||u_0||_{L^2(S^d)}||v_0||_{L^2(S^d)} 
$$

whenever $u_0$ is spectrally localized in the dyadic region $\sqrt{-\Delta_g} \in [N_1,2N_1)$, $v_0$ in the region $\sqrt{-\Delta_g} \in [N_2,2N_2)$, $N_2 \leq N_1$, with $\alpha=\frac{1}{4}+\epsilon$ when $d=2$ and $\alpha=\frac{1}{2}+\epsilon$ when $d=3$.

Using Theorem \ref{Bilinear FIO estimate}, we will be able to prove the following bilinear estimate for any closed manifold $(M,g)$:

\begin{theorem}\label{BLS on M theorem} Suppose $u_0, v_0\in L^2(M^d)$ are spectrally localized at dyadic scales $N_1$ and $N_2$ as above with $N_2 \leq N_1$. Then the following estimate holds:

\begin{equation}\label{BLS on M}
\left|\left|e^{it\Delta_g}u_0 e^{it\Delta_g}v_0\right|\right|_{L^2_{t,x}([-\frac{1}{N_1},\frac{1}{N_1}]\times M)} \lesssim_M \frac{N_2^{\frac{d-1}{2}}}{N_1^{\frac{1}{2}}}||u_0||_{L^2(M)}||v_0||_{L^2(M)}
\end{equation}

More generally,

\begin{equation}\label{BLS on M T}
||e^{it\Delta_g} u_0 e^{it\Delta_g}v_0||_{L^2([-T,T]\times M)} \leq \Lambda(T,N_1,N_2)||u_0||_{L^2(M)}||v_0||_{L^2(M)}
\end{equation}

where

\begin{equation}\label{Lambda T}
\Lambda(T,N_1,N_2) \lesssim_M \left\{
    \begin{array}{ll}
        \frac{N_2^{\frac{d-1}{2}}}{N_1^{\frac{1}{2}}}  & \mbox{if } T \ll N_1^{-1} \\
       T^{\frac{1}{2}}N_2^{\frac{d-1}{2}} & \mbox{if } T \gtrsim N_1^{-1}
    \end{array}
\right.
\end{equation}

In particular, for $T=1$ we have:

\begin{equation}\label{BLS on M T=1}
||e^{it\Delta_g} u_0 e^{it\Delta_g}v_0||_{L^2([-1,1]\times M)} \lesssim N_2^{(d-1)/2}||u_0||_{L^2(M)}||v_0||_{L^2(M)}
\end{equation}

\end{theorem}

Some notes are in order: First we notice that in the semiclassical/ short-range case $\eqref{BLS on M}$, the coefficient $\frac{N_2^{\frac{d-1}{2}}}{N_1^{\frac{1}{2}}}$ is the same as that on $\R^d$. This conforms with the heuristic that in the time interval $[0,\frac{1}{N_1}]$, the two waves $e^{it\Delta_g}v_0$ (which is moving with speed $\sim N_1$) and $e^{it\Delta_g}v_0$ (moving at speed $\sim N_2\leq N_1$) do not leave a coordinate patch and hence their product satisfies the same estimate as that on $\R^d$. Second, the estimates in $\eqref{BLS on M T}$ and $\eqref{BLS on M T=1}$ are essentially obtained from $\eqref{BLS on M}$ by splitting the time interval into pieces of length $N_1^{-1}$. It should be emphasized though that the exact dependence of $\Lambda(T, N_1, N_2)$ on its all parameters is often of great importance in applications (see \cite{H}). In fact, it is easy to see that bilinear estimates on the interval $[0,T]$ translate by scaling into bilinear estimates on the interval $[0,1]$ for the rescaled manifold $\lambda M$\footnote{Here $\lambda M$ can either be viewed as the Riemmannian manifold $(M, \frac{1}{\lambda^2}g)$ or by embedding $M$ into some ambient space $R^N$ and then applying a dilation by $\lambda$ to get $\lambda M$.}. The $\lambda-$dependence of those estimates is dictated by dependence of $\Lambda(T,N_1,N_2)$ on all its parameters. The bilinear Strichartz estimates on $\lambda M$ take the following form (see \cite{H} for relevant calculations):

\begin{corollary}\label{Strichartz on rescaled manifold}(Time $T$ estimate on $M$ implies time $1$ estimate on $\lambda M$)

Let $M$ be a $2D$ closed manifold and suppose that $N_1,N_2 \in 2^{\Z}$ and suppose $u_0,v_0 \in L^2(\lambda M)$
are spectrally localized around $N_1$ and $N_2$ respectively, with
$N_2\leq N_1$. Then

\begin{eqnarray}\label{bilinear estimate on lambda M}
||e^{it\Delta_{\lambda}} u_0 e^{it\Delta_{\lambda}}v_0||_{L^2([0,1]\times \lambda M)}
\lesssim_M& \Lambda(\lambda^{-2},\lambda N_1,\lambda N_2)||u_0||_{L^2(\lambda M)}||v_0||_{L^2(\lambda M)}\\
\lesssim_M& \left\{
    \begin{array}{ll}
        \left(\frac{N_2}{N_1}\right)^{1/2} ||u_0||_{L^2(\lambda M)}||v_0||_{L^2(\lambda M)} & \mbox{if } \lambda \gg N_1 \\
       \left(\frac{N_2}{\lambda}\right)^{1/2}||u_0||_{L^2(\lambda M)}||v_0||_{L^2(\lambda M)} & \mbox{if } \lambda \lesssim N_1
    \end{array}
\right.\label{the N_2/lambda decay}
\end{eqnarray}

where we have denoted by $\Delta_\lambda$ the Laplace-Beltrami operator on the rescaled manifold $\lambda M$.
\end{corollary}

Having favorable bounds (in terms of $\lambda$ and $N_2$) on the right hand side of $\eqref{the N_2/lambda decay}$ is crucial to obtaining global well-posedness of some nonlinear equations on $M$ below energy norm. In fact, in \cite{H} it is proven that the cubic nonlinear Schr\"odinger equation is globally well-posed in $H^s(M)$ for any closed $2D$ surface $M^2$ and all $s>2/3$, a result which matches the current (to the best of our knowledge) minimum regularity needed for global well-posedness on the 2-torus.

Finally, we note that as in the case of bilinear estimates on $\R^n$, the bilinear estimates in $\eqref{BLS on M}$ and $\eqref{BLS on M T}$ offer a refinement to those obtained by using linear estimates alone. However, this refinement is only visible when one looks at estimates over time intervals $[0,T]$ for $T\ll N_2^{-1}$ (or alternatively, estimates on rescaled manifolds). For example, for $d \geq 3$, applying H\"older's inequality, the $L_t^\infty L_x^2$ bound on $e^{it\Delta}u_0$, Bernstein and the $L_t^2L_x^{\frac{2d}{d-2}}$ for $e^{it\Delta}v_0$, one gets:

$$
||e^{it\Delta}u_0 e^{it\Delta}v_0||_{L_{t,x}^2([0,T]\times M)}\lesssim C(T, N_2)||u_0||_{L^2(M)}||v_0||_{L^2(M)}
$$

where $C(T, N_2)=N_2^{\frac{d-2}{2}}=\frac{N_2^{\frac{d-1}{2}}}{N_2^{1/2}}$ for $T\lesssim N_2^{-1}$ and $C(T,N_2)=T^{1/2}N_2^{\frac{d-1}{2}}$ for $T \geq N_2^{-1}$. This shows the improvement offered by $\eqref{BLS on M T}$ in the range $T \ll N_2^{-1}$ (especially when dealing with low-high frequency interaction $N_2 \ll N_1$). This improvement is due to the cancellation happening when we multiply the high frequency wave with the low frequency one. This cancellation is completely ignored by linear estimates. In the case, $d=2$, one would need to prove an estimate for the inadmissible pair $(q,r)=(2,\infty)$. This is possible with an $N^\epsilon$ loss. See \cite{jiang}. In this case, the bilinear estimate $\eqref{BLS on M T}$ not only offers a refinement to linear estimates at time scales $T \ll 1$  and in the range $N_2 \ll N_1$, but also yields better estimates in the time scale $T=1$ (no $N_2^\epsilon$ loss in $\eqref{BLS on M T=1}$). See \cite{H} for details.

The paper is organized as follows. In section \ref{proof of main theorem} we provide the proof of Theorem \ref{Bilinear FIO estimate}. In section \ref{proof of BLS on M}, we review the needed facts about the parametrix construction in \cite{BGT} and prove Theorem \ref{BLS on M theorem}. Finally in section \ref{further results} we prove inhomogeneous versions of the bilinear Strichartz estimates stated above in addition to mixed type bilinear estimates for products of the Schr\"odinger propagator $e^{it\Delta}u_0$ and the half wave propagators $e^{\pm it|\nabla|v_0}$. These estimates can also be deduced from Theorem \ref{Bilinear FIO estimate} and have potential applications (to be investigated elsewhere) in studying Zakharov type systems on closed manifolds. We use the notation $A \lesssim B$ to denote $A \leq C B$ for some $C>0$ and $A\sim B$ to denote $A\lesssim B \lesssim C$.

\textbf{Acknowledgements:} The author is deeply grateful to his advisor, Prof. Terence Tao, for his invaluable support, encouragement, and guidance. He also wishes to extend his immense gratitude to the referee for his careful review of the manuscript and his helpful comments and suggestions that considerably improved and clarified the exposition.

\section{Proof of Theorem \ref{Bilinear FIO estimate}}\label{proof of main theorem}

All implicit constants are allowed to depend on $d, \delta$ and uniform bounds on a finite number of derivatives of $\phi, \psi, a,$ and $b$.

\begin{equation}\label{product 1}
T_\lambda f(t,x)\tilde T_\mu g(t,x)=\int_{\R^d}\int_{\R^d}e^{i\left(\lambda \phi(t,x,\xi_1)+\mu \psi(t,x,\xi_2)\right)}a(t,x,\xi_1)b(t,x,\xi_2)f(\xi_1)g(\xi_2)d\xi_1d\xi_2.
\end{equation}

Since the supports of $a$ and $b$ are compact, one can use a finite partition of unity to split $a$ and $b$ into finitely many pieces so that on the support of each piece there exists $t_0,x_0,\xi_0,\xi_{2,0}$ such that 

$$
|t-t_0|, |x-x_0|,|\xi_1-\xi_0|,|\xi_2-\xi_{2,0}|\leq \frac{1}{C}
$$

where $C$ is some large constant depending only on $\delta$ and the uniform norms of $\phi$ and $\psi$ and their derivatives on the compact supports of $a$ and $b$.

Also notice that by applying a rotation $L$ of the domain $\R\times \R^d$: $(t,x)=L^{T}(s,y)$, the left hand side of $\eqref{bilinear estimate}$ is unaffected, whereas the hyper-surfaces $S_\phi$ and $S_\psi$ are both rotated by $L$. In fact, since:

$$
\nabla_{(s,y)}\left(\phi(L^{T}(s,y),x,\xi)\right)=L\left(\nabla \phi\right)(L^{T}(s,y),\xi)
$$

where $\nabla$ is taken in the first $n+1$ variable of $\phi$. Consequently, if we apply the change of variable $(t,x)=L^T(s,y)$, the canonical hyper-surfaces $S_\phi$ and $S_\psi$ are both rotated by $L$. Using this symmetry, one can assume that on the support of $a$ (resp. $b$):

\begin{equation}\label{rotated nondegeneracy}
\left|\det\left(\frac{\partial^2 \phi}{\partial \xi \partial x}(t_0,x_0,\xi_0)\right)\right|, \left|\det\left(\frac{\partial^2 \psi}{\partial \xi \partial x}(t_0,x_0,\xi_{2,0})\right)\right| \gtrsim 1.
\end{equation}

This means that the surfaces $S_\phi$ and $S_\psi$ can be regarded as graphs of functions of the form $(\xi,\tau_1(\xi))$ and $ (\xi,\tau_2(\xi))\subset T^*_{(t_0,x_0)}\R^{n+1}$ respectively.

Define:

$$
A:=\frac{\partial^2 \phi}{\partial \xi \partial x}(t_0,x_0,\xi_0); \;\;\;\;B:=\frac{\partial^2 \psi}{\partial \xi \partial x}(t_0,x_0,\xi_{2,0}).
$$

By the above, we have that $A$ and $B$ are invertible. It will be convenient later on to do the following change of variables in the $\xi_1$ integral and define $\xi=\xi_1+\frac{\mu}{\lambda}A^{-1}B\xi_2$\footnote{The justification for this change of variables will be obvious later on. However, at a heuristic level this corresponds to adding the momenta of the two waves.}. This gives:
\begin{equation}\label{product 2}
T_\lambda f(t,x)\tilde T_\mu g(t,x)=\int_{\R^d}\int_{\R^d}e^{i\lambda\left( \phi(t,x,\xi-\frac{\mu}{\lambda}A^{-1}B\xi_2)+\frac{\mu}{\lambda} \psi(t,x,\xi_2)\right)}c(t,x,\xi,\xi_2)f(\xi-\frac{\mu}{\lambda}A^{-1}B\xi_2)g(\xi_2)d\xi d\xi_2
\end{equation}

where we denoted $c(t,x,\xi,\xi_2)=a(t,x,\xi-\frac{\mu}{\lambda}A^{-1}B\xi_2)b(t,x,\xi_2)$ and all we have to remember about $c$ is that it is uniformly bounded along with all its derivatives (since $\frac{\mu}{\lambda}\leq 1$) and is supported in a small neighborhood of $(t_0,x_0,\xi_0+\frac{\mu}{\lambda}A^{-1}B \xi_{2,0},\xi_{2,0})$ of diameter $\lesssim \frac{1}{C}$. In particular, we have:

\begin{equation}\label{small support near xi_0}
|\xi-\frac{\mu}{\lambda}A^{-1}B\xi_2-\xi_0|\leq \frac{1}{C}
\end{equation}

for every $\xi,\xi_2$ in the support of $c$.

We now fix a particular coordinate direction $e_j$ (to be specified later), and write $\xi_2=p e_j+ \xi_2'$. Roughly speaking, the direction will be chosen using the transversality assumption of the two surfaces $S_\phi$ and $S_\psi$ so that 

$$
|\langle \nu_1(\xi_0), \frac{\partial^2 \psi(t_0,x_0,\xi_{2,0})}{\partial \xi \partial (t,x)} e_j \rangle|\gtrsim_\delta 1.\footnote{$\nu_1(\xi_0)\in \R^{n+1}$, $\frac{\partial^2 \psi(t_0,x_0,\xi_{2,0})}{\partial \xi \partial (t,x)}$ is an $(n+1)\times n$ matrix, and $e_j \in \R^n$, so the above expression makes sense.
}
$$

This will be possible because $\nu_2$ is the unique direction for which $\langle \nu_2,\frac{\partial^2 \psi(t_0,x_0,\xi_{2,0})}{\partial \xi \partial (t,x)}\rangle=\vec{0}_{\R^d}$ and since $\nu_1$ is quantitatively distinct from $\nu_2$, the vector $\langle \nu_1,\frac{\partial^2 \psi(t_0,x_0,\xi_{2,0})}{\partial \xi \partial (t,x)}\rangle\neq \vec{0}_{\R^d}$ and hence there exists a coordinate direction $e_j$ onto which the projection of this nonzero vector does not vanish. In other words, $\langle\nu_1,\frac{\partial^2 \psi(t_0,x_0,\xi_{2,0})}{\partial \xi \partial (t,x)}e_j\rangle$ can be thought of as the projection of $\nu_1$ onto the curve in $S_\phi(t_0,x_0)$ given by $t\mapsto \nabla_{(t,x)}\psi(t_0,x_0,\xi_{2,0}+te_j)$.

For convenience of notation, when confusion does not arise, we will assume that $j=1$ and write $\xi_2=(p,\xi_2')$ where $p \in \R$ and $\xi_2' \in \R^{d-1}$. As a result, we have:

\begin{align*}
&\left|\left|T_\lambda f(t,x)\tilde T_\mu g(t,x)\right|\right|_2\\
&=\left|\left|\int_{\R^{d-1}_{\xi'}}\int_{\R^d_\xi}\int_{\R_p}e^{i\lambda\left( \phi(t,x,\xi-\frac{\mu}{\lambda}A^{-1}B\xi_2)+\frac{\mu}{\lambda} \psi(t,x,\xi_2)\right)}c(t,x,\xi,\xi_2)f(\xi-\frac{\mu}{\lambda}\xi_2)g(\xi_2)d\xi dp d\xi_2' \right|\right|_2\\
&\leq\int_{\R^{d-1}_{\xi_2'}}\left|\left|\int_{\R^d_\xi}\int_{\R_p}e^{i\lambda\left( \phi(t,x,\xi-\frac{\mu}{\lambda}A^{-1}B\xi_2)+\frac{\mu}{\lambda} \psi(t,x,\xi_2)\right)}c(t,x,\xi,\xi_2)f(\xi-\frac{\mu}{\lambda}\xi_2)g(\xi_2)d\xi dp\right|\right|_{L^2_{t,x}}d\xi_2'.
\end{align*}

Freezing $\xi_2'$, we define the operator $S=S_{\xi_2'}: L^2(\R^{d+1})\to L^2(\R^{d+1})$ given by:

\begin{equation}\label{def of S}
SF(t,x)=\int_{\R^d_\xi}\int_{\R_p}e^{i\lambda\left( \phi(t,x,\xi-\frac{\mu}{\lambda}A^{-1}B\xi_2)+\frac{\mu}{\lambda} \psi(t,x,\xi_2)\right)}c(t,x,\xi,\xi_2)F(\xi,p) d\xi dp
\end{equation}

where $\xi_2=(p,\xi_2')$. As a result of this definition, our estimate is reduced to proving that for each $\xi_2'$, the following estimate holds for $S$:

\begin{equation}\label{estimate on S}
||SF||_{L^2_{t,x}(\R^{d+1})} \lesssim \frac{1}{\lambda^{d/2}\mu^{1/2}}||F||_{L^2_{p,\xi}(\R^{d+1})}.
\end{equation}

In fact, with such an estimate and by Cauchy-Schwarz in the $\xi_2'$ integral (keeping in mind that $c$ is compactly supported), we get that:

\begin{align*}
\left|\left|T_\lambda f(t,x)\tilde T_\mu g(t,x)\right|\right|_2
\lesssim& \frac{1}{\lambda^{d/2}\mu^{1/2}} \int_{|\xi_2'|\lesssim 1}\left|\left| f(\xi-\frac{\mu}{\lambda}(p,\xi_2'))g(p,\xi_2')\right|\right|_{L^2_{p,\xi}}d\xi_2'\\
\lesssim& \frac{1}{\lambda^{d/2}\mu^{1/2}}||f||_{L^2}||g||_{L^2}.
\end{align*}

The bound on $S$ is proved using a $T^*T$ argument. For convenience of notation, let us define:

\begin{equation}\label{def of Phi}
\Phi(t,x,\xi,p)= \phi(t,x,\xi-\frac{\mu}{\lambda}A^{-1}B\xi_2)+\frac{\mu}{\lambda} \psi(t,x,\xi_2)
\end{equation}

where $\xi_2=(p,\xi_2')$. With this notation, $S$ takes the form:

$$
SF(t,x)=\int_{\R^d_\xi}\int_{\R_p}e^{i\lambda\Phi(t,x,\xi,p)}c(t,x,\xi,p)F(\xi,p) d\xi dp.
$$

The adjoint of $S$ is given by the operator:

\begin{equation*}
S^*G(\xi,p)=\int_{\R^d_x}\int_{\R_t} e^{-i\lambda\Phi(t,x,\xi,p)}
\bar c(t,x,\xi,p)G(x,t)dxdt.
\end{equation*}

As a result, we get that:

\begin{equation}\label{def of S^*S}
S^*SF(\zeta,q)=\int_{\R^d_\xi}\int_{\R_p} K(\zeta, q,\xi,p) F(\xi,p)d\xi dp
\end{equation}

where 

\begin{equation}\label{def of K}
K(\zeta,q,\xi,p)=\int_{\R_t}\int_{\R^d_x} e^{i\lambda[\Phi(t,x,\xi,p)-\Phi(t,x,\zeta,q)]}c(t,x,\xi,p)\bar c(t,x,\zeta,q) dx dt.
\end{equation}

Our aim will be to show that $K$ satisfies the following bound:

\begin{equation}\label{bound on K}
K(\zeta,q,\xi,p) \lesssim_N \frac{1}{\left(1+\lambda|\xi-\zeta|+\mu|q-p|\right)^N}
\end{equation}

for a sufficiently large $N$ (any $N> d+1$ would do).

In fact, with such an estimate, one can easily see (using Schur's test for example) that $||S^*S||_{L^2 \to L^2} \lesssim \frac{1}{\lambda^d\mu}$. Since $||S||_{L^2 \to L^2}=||S^*S||_{L^2\to L^2}^{1/2}$ one gets that $||S||_{L^2 \to L^2}$ is bounded by $O(\frac{1}{\lambda^{d/2}\mu^{1/2}})$.

The bound on $K$ is  based on non-stationary--phase--type estimates and integration by parts. These are based on the following estimates on the phase function $\Phi$ and its derivatives.

\begin{lemma}\label{estimates on nabla Phi}
There exists $\Omega \in S^d$ such that:

1)
\begin{equation}\label{lower bound on nabla Phi}
\left|\langle\nabla_{t,x}\Phi(t,x,\xi,p)-\nabla_{t,x}\Phi(t,x,\zeta,q),\Omega\rangle \right| \gtrsim |\xi-\zeta|+\frac{\mu}{\lambda}|p-q|.
\end{equation}

2)
\begin{equation}\label{upper bound on delta Phi}
\left|\frac{\partial}{\partial x^{\alpha}\partial t^{\beta}}\left(\Phi(t,x,\xi,p)-\Phi(t,x,\zeta,q)\right)\right| \lesssim_{\alpha,\beta} |\xi-\zeta|+\frac{\mu}{\lambda}|p-q|.
\end{equation}
\end{lemma}

\proof
The second estimate $\eqref{upper bound on delta Phi}$ is a direct consequence of the definition $\eqref{def of Phi}$, the Taylor expansion, and the uniform boundedness of all the the $t,x$ derivatives of $\phi$ and $\psi$. We now turn to the proof of $\eqref{lower bound on nabla Phi}$.

Here we split the analysis into two cases:

\subsection{\emph{Case 1: $|\xi-\zeta|\geq \frac{1}{100}\frac{\mu}{\lambda}|p-q|$}}

The change of variables we have made in $\eqref{product 2}$ will allow us to prove $\eqref{lower bound on nabla Phi}$ in this case using only the $x$ derivative part of $\nabla_{t,x} \Phi$. In fact, using $\eqref{def of Phi}$, we have:

\begin{align}
\nabla_{x}\Phi(t,x,\xi,p)-\nabla_{x}\Phi(t,x,\zeta,q)=&\nabla_{x}\phi(t,x,\xi-\frac{\mu}{\lambda}A^{-1}B\xi_2)-\nabla_{x}\phi(t,x,\zeta-\frac{\mu}{\lambda}A^{-1}B\zeta_2)\label{delta phi}\\
&+\frac{\mu}{\lambda}\left(\nabla_{x}\psi(t,x,\xi_2)-\nabla_{x}\psi(t,x,\zeta_2)\right)\label{delta psi},
\end{align}
where $\zeta_2=(q, \xi_2')$.
We estimate $\eqref{delta phi}$ in the following manner:

\begin{align*}
\nabla_x \phi(t,x,\xi-\frac{\mu}{\lambda}A^{-1}B\xi_2)-\nabla_x \phi(t,x,\zeta-\frac{\mu}{\lambda}A^{-1}B\zeta_2)=&
\langle \frac{\partial^2 \phi}{\partial \xi \partial x}(t,x,\xi-\frac{\mu}{\lambda}A^{-1}B\xi_2),\xi-\zeta-\frac{\mu}{\lambda}A^{-1}B(\xi_2-\zeta_2)\rangle\\ 
&+O(|\xi-\zeta|^2)\\
=&\langle \frac{\partial^2 \phi}{\partial \xi \partial x}(t_0,x_0,\xi_0),\xi-\zeta-\frac{\mu}{\lambda}A^{-1}B(\xi_2-\zeta_2)\rangle+ \operatorname{Error}_1\\
=&A(\xi-\zeta)-\frac{\mu}{\lambda}B(\xi_2-\zeta_2)+\operatorname{Error}_1.
\end{align*}

where we used the fact that $A=\frac{\partial^2 \phi}{\partial \xi \partial x}(t_0,x_0,\xi_0)$. Here the $\operatorname{Error}_1$ term is:

\begin{align*}
\operatorname{Error}_1=&\left\langle \frac{\partial^2 \phi}{\partial \xi \partial x}(t,x,\xi-\frac{\mu}{\lambda}A^{-1}B\xi_2),\xi-\zeta-\frac{\mu}{\lambda}A^{-1}B(\xi_2-\zeta_2)\right\rangle\\
&-\left\langle \frac{\partial^2 \phi}{\partial \xi \partial x}(t_0,x_0,\xi_0),\xi-\zeta-\frac{\mu}{\lambda}A^{-1}B(\xi_2-\zeta_2)\right\rangle+O(|\xi-\zeta|^2).\\
\end{align*}

By our assumption of smallness of the support of $c$ (cf. $\eqref{small support near xi_0}$), the error can be estimated (if $C$ is chosen large enough depending on the uniform norms of derivatives of $\phi$) by:

$$
|\operatorname{Error}_1|\lesssim_{\phi}\frac{1}{C}|\xi-\zeta-\frac{\mu}{\lambda}A^{-1}B(\zeta_2-\xi_2))|+O(|\xi-\zeta|^2)\leq \frac{\gamma_1}{10}|\xi-\zeta|
$$

where $\gamma_1$ is chosen to be the smallest singular value of $A$ (or equivalently $\gamma_1=\min_{z\in S^{d-1}}|A z|$).

Next we estimate $\eqref{delta psi}$:

\begin{align*}
\frac{\mu}{\lambda}\left(\nabla_{x}\psi(t,x,\xi_2)-\nabla_{x}\psi(t,x,\zeta_2)\right)&=\frac{\mu}{\lambda}\left\langle\frac{\partial^2 \psi}{\partial \xi \partial x}(t,x,\xi_2),\xi_2-\zeta_2\right\rangle+O(\frac{\mu}{\lambda}|\xi_2-\zeta_2|^2)\\
&=\frac{\mu}{\lambda}\left\langle\frac{\partial^2 \psi}{\partial \xi \partial x}(t_0,x_0,\xi_{2,0}),\xi_2-\zeta_2\right\rangle+\operatorname{Error}_2\\
&=\frac{\mu}{\lambda}B(\xi_2-\zeta_2)+\operatorname{Error}_2
\end{align*}

where 
\begin{align*}
\operatorname{Error}_2=&\frac{\mu}{\lambda}\left(\left\langle\frac{\partial^2 \psi}{\partial \xi \partial x}(t,x,\xi_2),\xi_2-\zeta_2\right\rangle-\left\langle\frac{\partial^2 \psi}{\partial \xi \partial x}(t_0,x_0,\xi_{2,0}),\xi_2-\zeta_2\right\rangle\right)\\
&+O(\frac{\mu}{\lambda}|\xi_2-\zeta_2|^2),
\end{align*}

which, as before, can be bounded (using the fact that $|\xi_2-\zeta_2|,|\xi_2-\xi_{2,0}|\lesssim 1/C$ and that $\frac{\mu}{\lambda}|\xi_2-\zeta_2|\leq 100 |\xi-\zeta|$) by:

$$
|\operatorname{Error}_2|\leq \frac{\gamma_1}{10}|\xi-\zeta|.
$$

Collecting the above estimates we get that:

\begin{equation}\label{delta nabla_x Phi}
\nabla_x \Phi(t,x,\xi,p)-\nabla_x \Phi(t,x,\zeta,q)=A(\zeta-\xi)+\operatorname{Error}
\end{equation}

where $\operatorname{Error}=\operatorname{Error}_1+\operatorname{Error}_2$ is bounded by $\frac{\gamma_1}{5}|\zeta-\xi|$. We now let $\omega \in S^{d-1}$ be equal to $A(\zeta-\xi)/|A(\zeta-\xi)|$. Since

$$
|\langle A(\zeta-\xi),\omega\rangle| =|A(\xi-\zeta)| \geq \gamma_1|\xi-\zeta|
$$

by the definition of $\gamma_1$, we get that:
$$
\left|\langle \nabla_x \Phi(t,x,\xi,p)-\nabla_x \Phi(t,x,\zeta,q),\omega\rangle\right|\gtrsim |\xi-\zeta|.
$$

As a result, by taking $\Omega\in S^d$ equal to $(\omega,0)$ we get:

\begin{equation}
\left|\langle\nabla_{t,x}\Phi(t,x,\xi,p)-\nabla_{t,x}\Phi(t,x,\zeta,q),\Omega\rangle \right| \gtrsim |\xi-\zeta| \gtrsim |\xi-\zeta|+\frac{\mu}{\lambda}|p-q|
\end{equation}

which is $\eqref{lower bound on nabla Phi}$ in \emph{Case 1}.
\newline

\subsection{\emph{Case 2: ($|\xi-\zeta|\leq \frac{1}{100}\frac{\mu}{\lambda}|p-q|$):}}

The analysis in this case is a bit more delicate as it is here that the transversality assumption is used. In this case, we will take $\Omega=\nu_1(\xi_0)$, the normal to the surface $\xi \mapsto \nabla_{t,x}\phi(t_0,x_0,\xi)$ at $\xi_0$. With this choice we have:

\begin{align}
\left\langle\nabla_{t,x}\Phi(t,x,\xi,p)-\nabla_{t,x}\Phi(t,x,\zeta,q),\Omega\right\rangle =&\left\langle\nabla_{t,x}\phi(t,x,\xi-\frac{\mu}{\lambda}A^{-1}B\xi_2)-\nabla_{t,x}\phi(t,x,\zeta-\frac{\mu}{\lambda}A^{-1}B\zeta_2),\nu_1(\xi_0)\right\rangle \label{first line case 2}\\
&+\frac{\mu}{\lambda}\left\langle \nabla_{t,x}\psi(t,x,\xi_2)-\nabla_{t,x}\psi(t,x,\zeta_2),\nu_1(\xi_0)\right\rangle.\label{second line case 2}
\end{align}

The main term in this expression comes from $\eqref{second line case 2}$, whereas $\eqref{first line case 2}$ will be treated as an error. We start by lower bounding $\eqref{second line case 2}$.

Since 

\begin{align*}
\nabla_{t,x}\psi(t,x,\xi_2)-\nabla_{t,x}\psi(t,x,\zeta_2)=&\left\langle\frac{\partial^2\psi}{\partial \xi \partial (x,t)}(t,x,\xi_2),\xi_2-\zeta_2\right\rangle +O(|\xi_2-\zeta_2|^2)\\
=&\left\langle\frac{\partial^2\psi}{\partial \xi \partial (x,t)}(t_0,x_0,\xi_{2,0}),\xi_2-\zeta_2\right\rangle+ \operatorname{Error}_1\\
=&(p-q)\left\langle\frac{\partial^2\psi}{\partial \xi \partial (x,t)}(t_0,x_0,\xi_{2,0}),e_j \right\rangle +\operatorname{Error}_1\\
\end{align*}

where 

$$
\operatorname{Error}_1=\left\langle\frac{\partial^2\psi}{\partial \xi \partial (x,t)}(t,x,\xi_2),\xi_2-\zeta_2\right\rangle -\left\langle\frac{\partial^2\psi}{\partial \xi \partial (x,t)}(t_0,x_0,\xi_{2,0}),\xi_2-\zeta_2\right\rangle+O(|\xi_2-\zeta_2|^2).
$$

This is estimated as before using the small support assumption to get:

\begin{equation}\label{case 2 error 1}
|\operatorname{Error}_1|\lesssim_\psi \frac{1}{C}|\xi_2-\zeta_2|\leq \frac{1}{C}|p-q|
\end{equation}

where we have used in the last inequality the fact that $\xi_2=(p,\xi_2')$ and $\zeta_2=(q,\xi_2')$. We remark that the matrix $\frac{\partial^2\psi}{\partial \xi \partial (x,t)}(t_0,x_0,\xi_{2,0})$ is an $(n+1)\times n$ matrix and hence $\langle\frac{\partial^2\psi}{\partial \xi \partial (x,t)}(t_0,x_0,\xi_{2,0}),e_j\rangle$ is a vector in $\R^{n+1}$. From a geometric point of view, this vector lies in the tangent space to $S_\psi(t_0,x_0)$ at $\xi_{2,0}$.

Let us denote the $(n+1)\times n$ matrix

$$
N:=\frac{\partial^2\psi}{\partial \xi \partial (x,t)}(t_0,x_0,\xi_{2,0}).
$$

Recall that by definition, $\nu_2:=\nu_2(\xi_{2,0})$ is the unique vector (up to sign) in $S^d$ such that  $\nu_2^TN=0$ where $\nu_2^T$ is the row vector corresponding to $\nu_2$. In particular, the map from the $n-$dimensional subspace $\nu_2^\perp \subset \R^{n+1}$ into $\R^n$ given by:

$$
\nu\in \nu_2^{\perp} \mapsto \nu^TN\in \R^n
$$

is an isomorphism. Let $\gamma_2>0$ denote its smallest singular value (or equivalently $\gamma_2$ is the \emph{positive} infimum of the above map when $\nu \in \nu_2^{\perp}$ satisfies $||\nu||=1$). 

Writing $\nu_1(\xi_0)=\alpha \nu_2+\beta \nu_3$ with $\nu_3 \in \nu_2^{\perp}$,$||\nu_3||=1$, and $|\alpha|,|\beta|\leq 1$, we notice that since $1-\delta> |\langle \nu_1,\nu_2\rangle|=|\alpha|$ we have that $|\beta|=\sqrt{1-\alpha^2}\geq \sqrt{\delta}$. 

As a result, we have:

\begin{align*}
\left\langle \nu_1,\nabla_{t,x}\psi(t,x,\xi_2)-\nabla_{t,x}\psi(t,x,\zeta_2)\right\rangle=&
(p-q)\nu_1^TNe_j  +\operatorname{Error_1}
=\beta(p-q)\nu_3^TNe_j  + \operatorname{Error}_1.\\
\end{align*}

Since $||\nu_3^TN||\geq \gamma_2$, one can choose $e_j$ so that $|\nu_3^TNe_j|\geq \gamma_2/\sqrt d=:c_1$. Combining this to the estimate on $\operatorname{Error}_1$ in $\eqref{case 2 error 1}$ above we get that if $C$ is large enough:

\begin{equation}\label{second line case 2 estimate}
\left|\langle \nu_1,\nabla_{t,x}\psi(t,x,\xi_2)-\nabla_{t,x}\psi(t,x,\zeta_2)\rangle\right|\geq c_1 \sqrt \delta |p-q|  - \frac{c_1 \sqrt \delta}{100}|p-q|\geq \frac{99}{100}c_1\sqrt\delta |p-q|.
\end{equation}

As mentioned before, we will treat $\eqref{first line case 2}$ as an error. Indeed,

\begin{align*}
\left\langle\nabla_{t,x}\phi(t,x,\xi-\frac{\mu}{\lambda}A^{-1}B\xi_2)-\nabla_{t,x}\phi(t,x,\zeta-\frac{\mu}{\lambda}A^{-1}B\zeta_2),\nu_1(\xi_0)\right\rangle \\
=\nu_1(\xi_0)^TD_{(d+1)\times d}(t,x,\xi-\frac{\mu}{\lambda}A^{-1}B\xi_2)[\xi-\zeta-\frac{\mu}{\lambda}A^{-1}B(\xi_2-\zeta_2)]+O(|\frac{\mu}{\lambda}(p-q)|^2)
\end{align*}

where we have denoted 
$$
D_{(d+1)\times d}(t,x,\eta)=\frac{\partial^2\phi}{\partial \xi \partial (x,t)}(t,x,\eta)
$$

and also used that $|\xi-\zeta|\leq \frac{\mu}{\lambda}|p-q|$ in this case. Since the derivatives of $D$ are uniformly bounded and because of the small support assumption $\eqref{small support near xi_0}$, we have:

$$
||D_{(d+1)\times d}(t,x,\xi-\frac{\mu}{\lambda}A^{-1}B\xi_2)-D_{(d+1)\times d}(t_0,x_0,\xi_0)||\lesssim \frac{1}{C}\leq \frac{c_1\sqrt \delta}{100(||A^{-1}B||+1)}
$$
if $C$ is large enough.

Using the fact that $\nu_1^TD_{(d+1)\times d}(t_0,x_0,\xi_0)=0$, we get that 

\begin{equation}\label{first line case 2 estimate}
\left|\left\langle\nabla_{t,x}\phi(t,x,\xi-\frac{\mu}{\lambda}A^{-1}B\xi_2)-\nabla_{t,x}\phi(t,x,\zeta-\frac{\mu}{\lambda}A^{-1}B\zeta_2),\nu_1(\xi_0)\right\rangle\right|\leq \frac{c_1\sqrt\delta }{50}\frac{\mu}{\lambda}|p-q|
\end{equation}

again using the small support assumption. 

Combining $\eqref{first line case 2 estimate}$ and $\eqref{second line case 2 estimate}$, we get 
$\eqref{lower bound on nabla Phi}$ for \emph{Case 2}.

\endproof

Now we are ready to perform the integration by parts needed to prove the estimate $\eqref{bound on K}$.  Recall that 
$$
K(\zeta,q,\xi,p)=\int_{\R_t}\int_{\R^d_x} e^{i\lambda[\Phi(t,x,\xi,p)-\Phi(t,x,\zeta,q)]}c(t,x,\xi,p)\bar c(t,x,\zeta,q) dx dt.
$$

Let $D_\Omega$ be the operator given by: 

\begin{equation}\label{def of operator D}
D_{\Omega}:=\frac{1}{i \lambda \left\langle \nabla_{t,x}\Phi(t,x,\xi,p)-\nabla_{t,x}\Phi(t,x,\zeta,q), \Omega \right\rangle}\langle \nabla_{(x,t)}, \Omega\rangle.
\end{equation}

Then 
$$
D_{\Omega} \left(e^{i\lambda(\Phi(t,x,\xi,\xi_2)-\Phi(t,x,\zeta,\zeta_2)}\right)=e^{i\lambda(\Phi(t,x,\xi,\xi_2)-\Phi(t,x,\zeta,\zeta_2))}.
$$

Noticing that the formal adjoint of $D_{\Omega}$ acting on $L^2$ is:

$$
D_{\Omega}^T=\langle \nabla_{(x,t)}, \Omega\rangle \frac{1}{\left(i \lambda\langle \nabla_{t,x}\bar \Phi(t,x,\xi,p)-\nabla_{t,x}\bar \Phi(t,x,\zeta,q), \Omega \rangle\right)}
$$

we get that:

\begin{align*}
K(\zeta,q,\xi,p)=&\int_{\R_t}\int_{\R^d_x} e^{i\lambda[\Phi(t,x,\xi,p)-\Phi(t,x,\zeta,q)]}c(t,x,\xi,p)\bar c(t,x,\zeta,q) dx dt\\
&=\int_{\R_t}\int_{\R^d_x} e^{i\lambda[\Phi(t,x,\xi,p)-\Phi(t,x,\zeta,q)]}\overline{\left(D_\Omega^T\right)^N\bar c(t,x,\xi,p)c(t,x,\zeta,q)} dx dt.
\end{align*}

Using the estimates in Lemma $\eqref{estimates on nabla Phi}$, it is easy to see that that

$$
\left(D_\Omega^T\right)^N\bar c(t,x,\xi,p) c(t,x,\zeta,q)\lesssim_N \frac{1}{ \left(\lambda|\xi-\zeta|+\mu|p-q|\right)^N}.
$$

When $\lambda|\xi-\zeta|+\mu|p-q| \leq 1$, we do not perform any integration by parts and estimate the $K$ integrand by $O(1)$ and hence $K$ by $O(1)$ as well. Otherwise we use the above decay. As a result, we get that:

$$
K(\xi,\xi_2,\zeta,\zeta_2)\lesssim_N \frac{1}{\left(1+\lambda |\xi-\zeta|+\mu|p-q|\right)^N}
$$
which finishes the proof.

\endproof

\remark
It is not hard to see that the estimate $\eqref{bilinear estimate}$ is sharp. In fact, by considering the restriction case and taking $\phi(t,x,\xi)=\psi(t,x,\xi)=x.\xi+t|\xi|^2$ with $a$ having its $\xi$ support in the region $|\xi| \geq 100$ and $b$ having its $\xi$ support near $|\xi|\leq 1$, one can can reduce the sharpness of $\eqref{bilinear estimate}$ to that of $\eqref{BLS on R^d}$ which is known to be sharp. In fact, this can be seen by first reducing to the case when $N_2=1$ (again using scaling) and taking $\widehat{u_0}$ to be the characteristic function of $[N_1, N_1 +N_1^{-1}]\times [-1,1]^{d-1}$ (hence $||u_0||_{L_x^2}\sim N_1^{-1/2}$); and $\widehat{v_0}$ to be the characteristic function of $[-1,1]^d$ (hence $||v_0||_{L_x^2}\sim 1$). By Plancherel's theorem in space and time, we get that L.H.S of $\eqref{BLS on R^d} \gtrsim ||\chi_{R_1}*\chi_{R_2}||_{L^2(\R^{d+1})}$ where $R_1=[N_1, N_1+N_1^{-1}]\times [0,1]^d$ and $R_2=[-1,1]^{d+1}$. A direct calculation now shows that $\chi_{R_1}*\chi_{R_2} \gtrsim \frac{1}{N_1}\chi_{R_3}$ where $R_3=[N_1+\frac{1}{4}, N_1+\frac{3}{4}]\times [-\frac{1}{2},\frac{1}{2}]^d$ and hence $||\chi_{R_1}*\chi_{R_2}||_{L^2(\R^{d+1})} \sim \frac{1}{N_1}$, which gives that L.H.S of $\eqref{BLS on R^d} \gtrsim \frac{1}{N_1^{1/2}}||u_0||_{L_x^2}||v_0||_{L_x^2}$.

\section{Bilinear Strichartz Estimates}\label{proof of BLS on M}

We will apply the result of the previous section to get bilinear Strichartz estimates for the free Schr\"odinger evolution on compact manifolds without boundary. These will be analogues in the variable coefficient case to the estimate $\eqref{BLS on R^d}$ on $\R^d$ with the Euclidean Laplacian which we recall here for convenience:

\begin{equation*}
||e^{it\Delta}u_0 e^{it\Delta}v_0||_{L^2(\R \times \R^d)} \lesssim \frac{N_{2}^{\frac{d-1}{2}}}{N_{1}^{\frac{1}{2}}}||u||_{L^2(\R^d)}||v||_{L^2(\R^d)}
\end{equation*}

where $u,v \in L^2(\R^d)$ are frequency localized on the dyadic annuli $\{\xi\in \R^d: |\xi|\in[N_1,2N_1]\}$ and $\{\xi\in \R^d: |\xi|\in[N_2,2N_2]\}$ respectively. 

By scaling time and space, one can easily see that this estimate is equivalent to the same one on the time interval $[0,\frac{1}{N_1}]$. On this time scale, the numerology in $\eqref{BLS on R^d}$ can be understood (heuristically at least) by a simple back-of-the-envelope calculation. Thinking of $e^{it\Delta}u_0$ as a ``bump function" localized in frequency at scale $N_1$ and \emph{initially} (at $t=0$) localized in space at scale $\frac{1}{N_1}$. The evolution moves this bump function at a speed $N_1$ thus expanding its support at this rate while keeping the $L^2$ norm conserved. Similarly, $e^{it\Delta}v_0$ could be thought of as a ``bump function" that is initially concentrated in space at scale $\sim \frac{1}{N_2}$ and moving (expanding) at speed $N_2$. A simple schematic diagram allows to estimate the space-time overlap of the two expanding ``bump functions" thus giving the estimate $\frac{N_2^{(d-1)/2}}{N_1^{1/2}}$ for the $L^2_{t,x}([0,N_1^{-1}]\times \R^d)$ of the product.

The goal of this section is to prove the analogue of $\eqref{BLS on R^d}$ for the linear evolution of the Schr\"odinger equation on a $C^\infty$ compact manifold $M$ without boundary. This was stated in Theorem \ref{BLS on M theorem}. All implicit constants are allowed to depend on $M$ and the uniform bounds of its metric functions (they are all finite since $M$ is compact). To fix notation, we consider two functions $u_0,v_0 \in C^{\infty}(M)$\footnote{The full result for $u_0, v_0 \in L^2(M)$ can be obtained in the end by a standard limiting argument.} such that $u_0=\varphi(\frac{\sqrt {-\Delta}}{N_1})u_0$ and $v_0=\varphi(\frac{\sqrt {-\Delta}}{N_2})v_0$ where $\varphi\in C_0^{\infty}(\R)$, and we would like to estimate the $L^2_{t,x}$ norm of the product $e^{it\Delta}u_0e^{it\Delta}v_0$. We assume further that $\varphi$ vanishes in a small neighborhood of the origin.

\remark{The same analysis allows to consider different frequency localizations for $u_0$ and $v_0$ like $u_0=\varphi(\frac{\sqrt{-\Delta}}{N_1})u_0$ and $v_0=\psi(\frac{\sqrt{-\Delta}}{N_2})v_0$  with $\varphi,\psi\in C_0^{\infty}$ as long as $\varphi$ vanishes in a neighborhood of the origin and $N_1$ is sufficiently larger than $N_2$. In particular, $\psi$ does not need to vanish near the origin.} 

To simplify notation, we use $\Delta$ to denote the Laplace-Beltrami operator $\Delta_g$ on $M$,  and $|\xi|_{g(x)}$ to denote $\sqrt{g(x)^{ij}\xi_i\xi_j}$.

\emph{Proof of Theorem \ref{BLS on M theorem}}: The proof is organized as follows. We will first review some important facts about microlocalizing $\varphi(h\sqrt{-\Delta})$ and constructing the Schr\"odinger parametrix (as in \cite{BGT}) that will be used to approximate the linear evolutions. The case when $N_2\sim N_1$, will then follow directly from the semiclassical linear Strichartz estimates already proven in \cite{BGT}(Proposition 2.9). As a result, we will only need to consider the case when $N_2 \ll N_1$. This will ensure that the canonical hyper-surfaces associated to the phase functions of the parametrices are transversal as defined in the previous section, a fact which will allow us to apply Theorem \ref{Bilinear FIO estimate}.

\subsection{Microlocalizing $\varphi(h\sqrt{-\Delta})$\cite{BGT},\cite{Sogge},\cite{H2}}

In this section, we will briefly review how spectrally localizing a function $f\in C^{\infty}(M)$ using the spectral multiplier $\varphi(h\sqrt {-\Delta})$ is expressed in local coordinates. Essentially, up to  smooth remainder terms, $\varphi(h\sqrt{-\Delta})f$ is given in local coordinates as a pseudo-differential operator whose symbol $a(x,\xi)$ has a support that reflects the spectral localization dictated by $\varphi$:

\begin{proposition}\label{microlocalizing}
Let $\varphi \in C_0^\infty(\R)$ and $\kappa: U\subset \R^d \to V \subset M$ be a coordinate parametrization of $M$. Also let $\chi_1, \chi_2 \in C_0^\infty(V)$ be such that $\chi_2=1$ near the support of $\chi_1$. Then for every $N\in \N$, every $h\in(0,1)$, and every $\sigma \in [0,N]$, there exists $a_N(x,\xi)$ supported in $\{(x,\xi)\in U \times \R^d: \kappa(x)\in \operatorname{supp}(\chi_1), |\xi|_{g(x)} \in \operatorname{supp}(\varphi)\}$ such that:
\begin{equation}\label{microlocalizing equation}
\left|\left|\kappa^*\left(\chi_1\varphi(h\sqrt {-\Delta})f\right)-a(x,hD)\kappa^*(\chi_2f)\right|\right|_{H^\sigma(\R^d)}\lesssim_N h^{N-\sigma}||f||_{L^2(M)}
\end{equation}
for every $f\in C^\infty(M)$. In particular, if $\varphi$ is supported away from the origin, then so is the $\xi$ support of $a(x,\xi)$. Here $\kappa^*$ is used to denote the pull-back map given by: $\kappa^* f=f\circ \kappa$.
\end{proposition}

\proof See Proposition 2.1 of \cite{BGT} (alternatively, one can use the parametrix expression of the half-wave operator $e^{it\sqrt{-\Delta}}$ (see \cite{Sogge} for example), along with the expression of $\varphi$ in terms of its Fourier transform).

A consequence of this proposition and a finite partition of unity in $M$, one can split $u_0=\varphi(h\sqrt {-\Delta})u_0$ into pieces of the form $\chi_1 \varphi(h\sqrt {-\Delta})u_0$ and replace each of  those pieces (incurring an error that is $O(h^N||u_0||_{L^2})$) by $a(x,hD)\kappa^*(\chi_2 u_0)$ which is a compactly supported function in space and is pseudo-localized in frequency in the following sense:

\emph{There exists a function $\psi \in C_0^\infty(\R^d)$ such that for all $h\in (0,1),\sigma>0, \operatorname{and} N>0$,} 

\begin{equation}\label{microlocalization lemma}
\kappa^*(\chi_1 \varphi(h\sqrt {-\Delta})f)=\psi(hD)\kappa^*(\chi_1 \varphi(h\sqrt {-\Delta})f) +r_1
\end{equation} 

\emph{with $||r_1||_{H^\sigma(\R^d)}\lesssim _{\sigma, N}h^N||f||_{L^2}$. If $\varphi$ is supported away from 0, one can also take $\psi$ to be supported at a positive distance from the origin in $\R^d$.} This follows easily from Proposition \ref{microlocalizing} and standard pseudo-differential calculus (See for e.g. \cite{Stein}). We will denote $w_0(x)=a(x,hD)\kappa^*(\chi_2 u_0)$. In brief, $w_0$ is compactly supported in space and can be replaced by $\psi(hD)w_0$ at the cost of an error that is $O(h^N||u_0||_{L^2(M)})$.

\subsection{The Parametrix \cite{BGT}}
With this microlocalization setup, Burq, Gerard, and Tzvetkov constructed an approximate solution to the semiclassical equation:

\begin{eqnarray}\label{semiclassical equation}
ih\partial_t w +h^2\Delta_g w=&0\\
w(0)=&\varphi(h\sqrt{- \Delta})v_0
\end{eqnarray}

in local coordinates. More precisely, using the usual WKB construction (see for example \cite{H2},\cite{BGT}, or the lecture notes \cite{EZnotes}), they show that there exists $\alpha>0$, such that on the time interval $[-\alpha,\alpha]$ 

$$
w(s)=\tilde w(s) +r_2(s)
$$

where $r_2(s)$ satisfies $||r_2(t)||_{L_t^\infty([-\alpha,\alpha]\times H^\sigma(M))}\lesssim h^N||w_0||_{L^2(M)}$ (with $N$ sufficiently large) and $\tilde w (t)$ is supported in a compact subset of $V\subset M$ and is given in local coordinates by the following oscillatory integral:

\begin{equation}\label{main term of the parametrix}
\tilde{w}(s,x)=\frac{1}{(2\pi h)^d}\int_{\R^d}e^{\frac{i}{h}\tilde\phi(s,x,\xi)}a(s,x,\xi,h)\widehat w_0(\frac{\xi}{h})d\xi.
\end{equation}

Here $a(s,x,\xi,h)=\sum_{j=0}^Nh^ja_j(s,x,\xi)$, and $ a_j \in C_0^{\infty}([-\alpha,\alpha]\times U \times U'\subset\subset \R\times \R^d \times \R^d)$. $w_0$ is the microlocalization of $\varphi(h\sqrt \Delta)v_0$ described above. Since $w_0$ can be replaced by $\psi(hD)w_0$ at the cost of an error that is $O(h^N||w_0||_{L^2(\R^d)})$ one can assume without loss of generality that $a(s,x,\xi,h)$ has its $\xi$ support at a positive distance from the origin in frequency space if $\varphi$ is supported away from 0 itself.

The phase function $\tilde \phi$ appearing in the integral $\eqref{main term of the parametrix}$ satisfies the eikonal equation:

\begin{eqnarray}\label{eikonal equation}
\partial_s \tilde \phi +\sum_{ij}g^{ij}\partial_i \tilde \phi\partial_j \tilde \phi&=0\\
\tilde \phi(0,x,\xi)&=x.\xi.
\end{eqnarray}

\subsection{Semiclassical Linear Strichartz estimates and the case $N_1\sim N_2$} Using this representation, one can easily use stationary phase (see \cite{BGT} for details) to get the following semiclassical dispersion estimate:

\begin{equation}\label{semiclassical dispersion estimate}
||e^{it\Delta}\varphi^2(h\sqrt{ -\Delta})v_0||_{L^{\infty}(M)}\lesssim_M \frac{1}{t^{d/2}}||v_0||_{L^1(M)}
\end{equation}

for every $t\in [-\alpha h,\alpha h]$ with $0<\alpha \ll 1$. Combining this with the Keel-Tao machinery (see \cite{KT}) one immediately gets the following semiclassical Strichartz estimate:

\begin{equation}\label{semiclassical Strichartz estimate}
||e^{it\Delta}\varphi(h\sqrt {-\Delta})u_0||_{L_t^qL_x^r([-\alpha h,\alpha h]\times M)}\lesssim_M ||u_0||_{L^2(M)}
\end{equation}

whenever $2\leq q,r \leq \infty$ satisfy $\frac{2}{q}+\frac{d}{r}=\frac{d}{2}$ and $(q,r,d)\neq (2,\infty,2)$.

This estimate is enough to prove $\eqref{BLS on M}$ in the case when $h=\frac{1}{N_1}\sim m=\frac{1}{N_2}$. In fact, for $d=2$, one can use the $L_{t,x}^4$ Strichartz estimate to get:

$$
||e^{it\Delta}u_0 e^{it\Delta}v_0||_{L_{t,x}^2([-\alpha h,\alpha h]\times M^2)}\leq ||e^{it\Delta}\varphi(h\sqrt{- \Delta})u_0||_{L^4_{t,x}}||e^{it\Delta}\varphi(h\sqrt{- \Delta})v||_{L^4_{t,x}}\lesssim ||u_0||_{L^2(M^2)}||v_0||_{L^2(M^2)}.
$$

Whereas for $d \geq 3$, one can apply H\"older's inequality, the $L_t^\infty L_x^2$ bound on $e^{it\Delta}u_0$, Bernstein\footnote{One can verify Bernstein's inequality in the setting of compact manifolds by using Proposition 3.2 and the fact that the kernel $K(x,y)$ of $a(x,hD)$ satisfies the bound $\|K(x,y)\|_{L_x^r L_y^p(\R^d \times \R^d)}\lesssim_a h^{-d(1-\frac{1}{r}-\frac{1}{p})}$.} and the $L_t^2L_x^{\frac{2d}{d-2}}$ for $e^{it\Delta}v_0$ to get:

$$
||e^{it\Delta}u_0 e^{it\Delta}v_0||_{L_{t,x}^2([0,\alpha h]\times M)}\lesssim N_2^{\frac{d-2}{2}}||u_0||_{L^2(M)}||v_0||_{L^2(M)}
$$

as desired.

\subsection{The case $N_1 \gg N_2$}

In this section, we will reduce the case $N_1\gg N_2$ into a verification of the conditions of Theorem \ref{bilinear estimate}.
By rescaling time, we have:

\begin{equation}\label{rescaling time}
\begin{split}
||e^{it\Delta}u_0e^{it\Delta}v_0||_{L^2_{t,x}([-\alpha h,\alpha h]\times M)}= h^{1/2}||e^{iht\Delta}u_0e^{iht\Delta}v_0||_{L^2_{t,x}([-\alpha,\alpha]\times M)}\\
=h^{1/2}||e^{iht\Delta}u_0e^{im(\frac{h}{m}t)\Delta}v_0||_{L^2_{t,x}([-\alpha,\alpha]\times M)}.
\end{split}
\end{equation}

As a result it is enough to show that:
\begin{equation}\label{reduced Strichartz estimate}
||e^{iht\Delta}u_0e^{im(\frac{h}{m}t)\Delta}v_0||_{L^2_{t,x}([-\alpha,\alpha]\times M)}\lesssim \frac{1}{m^{\frac{d-1}{2}}}||u_0||_{L^2(M)}||v_0||_{L^2(M)}.
\end{equation}

The advantage of writing the estimate in this way is that we can now use the parametrices for $e^{ith\Delta}u_0$ and $e^{itm\Delta} v_0$ constructed above to write\footnote{Strictly speaking this representation only holds in an open neighborhood of $x_0\in M$. Since $M$ is compact, we can cover it by finitely many of such neighborhood, and hence we only need to prove the estimate on each one of them.}:

$$
e^{ith\Delta}u_0(x)=\tilde T_hu_0(t,x)+R_hu_0(t,x)
$$

and 
$$
e^{im(\frac{ht}{m})\Delta}v_0(x)=\tilde S_m v_0(t,x)+R_m v_0(t,x)
$$

where $\tilde T_h$ and $\tilde S_m$ are defined according to $\eqref{main term of the parametrix}$ by:

\begin{equation}\label{def of tilde T_h}
\tilde T_h u_0(t,x)=\frac{1}{(2\pi h)^d}\int_{\R^d}e^{\frac{i}{h}\tilde\phi(t,x,\xi)}a_1(t,x,\xi,h)\widehat{\tilde u_0}(\frac{\xi}{h})d\xi
\end{equation}

and 

\begin{equation}\label{def of tilde S_m}
\tilde S_m v_0(t,x)=\frac{1}{(2\pi m)^d}\int_{\R^d}e^{\frac{i}{m}\tilde\phi(\frac{ht}{m},x,\xi_2)}a_2(\frac{h}{m}t,x,\xi_2,m)\widehat {\tilde v_0}(\frac{\xi_2}{m})d\xi_2
\end{equation}

where $\tilde u_0$ and $\tilde v_0$ are the respective microlocalizations of $u_0$ and $v_0$ in the considered coordinate patch (in particular $||\tilde u_0||_{L^2(\R^d)}\lesssim ||u_0||_{L^2(M)}$ and $||\tilde v_0||_{L^2(M)}\lesssim ||v_0||_{L^2(M)}$). Also we have that:

\begin{equation}\label{Remainder terms}
||R_h u_0||_{L^\infty_t H^\sigma([-\alpha,\alpha]\times M)}\lesssim h^N||u_0||_{L^2(M)} \textrm{ and }
||R_m v_0||_{L^\infty_t H^\sigma([-\alpha,\alpha]\times M)}\lesssim m^N||v_0||_{L^2(M)}.
\end{equation}

The main contribution comes of course from the product $\tilde T_h u_0 \tilde S v_0$. For example the cross terms $\tilde T_h u_0 R_m v_0$ and $R_h u_0 \tilde S_m v_0$ can be bounded as follows:

$$
||\tilde T_h u_0 R_m v_0||_{L^2_{t,x}}\leq ||\tilde T_h u_0||_{L_t^\infty L_x^2}||R_m v_0||_{L_t^2 L_x^{\infty}}\lesssim ||u_0||_{L^2}||v_0||_{L^2}$$

where in the last step we used $\eqref{Remainder terms}$ and a crude Sobolev embedding to bound $||R_m v_0||_{L_t^2L_x^\infty}$ by $||R_m||_{L_t^2 H^\sigma_x}$ for some $\sigma> d/2$. The $L_t^\infty L_x^2$ bound on $\tilde T_h u_0$ follows from the the $L_t^\infty L_x^2$ boundedness of $e^{ith\Delta}u_0$. Similarly, one bounds the contributions of $R_h u_0 \tilde S_m v_0$ and $R_h u_0 R_m v_0$.

To bound the contribution of $\tilde T_h u_0 \tilde S_m v_0$, we now apply Theorem \ref{Bilinear FIO estimate} with $\phi(t,x,\xi)=\tilde \phi(t,x,\xi)$ and $\psi(t,x,\xi_2)=\tilde \phi(\frac{h}{m}t,x,\xi_2)$, $f(\xi):=\tilde u(\xi/h)$, and $g(\xi)=\tilde v_0(\xi/m)$, to get that:

$$
||\tilde T_h u_0 \tilde S_m v_0||_{L^2_{t,x}([-\alpha,\alpha]\times \R^d)}\lesssim \frac{1}{(hm)^d}(h^dm)^{1/2}||f||_{L^2(\R^d)}||g||_{L^2(\R^d)}\lesssim \frac{1}{m^{(d-1)/2}}||\tilde u_0||_{L^2(\R^d)}||\tilde v_0||_{L^2(\R^d)}
$$

which clearly gives $\eqref{reduced Strichartz estimate}$ and hence $\eqref{BLS on M}$. As a result, all we need to do is to verify that the requirements of Theorem $\eqref{Bilinear FIO estimate}$ are satisfied.

Obviously all derivatives of $\phi$ and $\psi$ are uniformly bounded on the compact supports of $a_1$ and $a_2$ ($\frac{h}{m}\leq 1$). Moreover, since $\tilde\phi(0,x,\xi)=x.\xi$, we have that $\frac{\partial^2 \phi}{\partial \xi \partial x}(0,x,\xi)=Id$ (invertible), the non-degeneracy condition $\eqref{nondegeneracy}$ is satisfied at $t=0$ and hence for all  $t\in [-\alpha,\alpha]$ if $\alpha$ is small enough. 

Now we consider the canonical surfaces $S_\phi$ and $S_\psi$:

Recall that $S_\phi$ and $S_\psi$ are the images of the maps:
\begin{align*}
\xi_1\mapsto \nabla_{t,x}\phi(t,x,\xi_1)&=(\nabla_x \tilde \phi(t,x,\xi_1), \partial_t \tilde \phi(t,x,\xi_1)) \\
\xi_2 \mapsto \nabla_{t,x}\psi(t,x,\xi_2)&=(\nabla_x \tilde \phi(\frac{h}{m}t,x,\xi_2), \frac{h}{m}\partial_t \tilde \phi(\frac{h}{m}t,x,\xi_2))
\end{align*}

respectively. By the non-degeneracy condition above, $S_\phi$ and $S_\psi$ are smooth embedded hyper-surfaces in $T^*_{(t,x)}\R^{n+1}$. We need to show that if $\nu_1(\xi_1)$ is the normal to $S_\phi$ at $\nabla_{t,x}\phi(t,x,\xi_1)$ and $\nu(\xi_2)$ is the normal to $S_\psi$ at $\nabla_{t,x}\psi(t,x,\xi_2)$, then there is a $\delta>0$ (uniform in $\xi_1$ and $\xi_2$) such that:

\begin{equation}\label{transversality schrodinger}
|\langle \nu_1, \nu_2 \rangle| \leq 1-\delta.
\end{equation}

By continuity, we only need to verify $\eqref{transversality schrodinger}$ at $t=0$ for all $x,\xi_1,\xi_2$. This will imply that the same holds for all $t\in [-\alpha,\alpha]$ if $\alpha$ is small enough. We now fix $(0,x_0)\in \R^{d+1}$ and consider the surfaces $S_\phi$ and $S_\phi$ in $T^*_{(0,x_0)}\R^{d+1}$. From the eikonal equation $\eqref{eikonal equation}$, $\tilde \phi(0,x,\xi)=x.\xi$ and $\partial_t \tilde \phi(0,x,\xi)=g^{ij}(x)\xi_i\xi_j$. A straight-forward computation gives that: 

$$
\nu_1(\xi)=\frac{(2g^{1j}\xi_j,2g^{2j}\xi_j,...,2g^{dj}\xi_j,-1)}{\sqrt{1+4|\xi|_{g(x)}^2}}
$$

and 

$$
\nu_2(\xi)=\frac{(2\frac{h}{m}g^{1j}\xi_j,2\frac{h}{m}g^{2j}\xi_j,...,2\frac{h}{m}g^{dj}\xi_j,-1)}{\sqrt{1+4|\frac{h}{m}\xi|_{g(x)}^2}},
$$

where we recall our notation that $|\xi|_{g(x)}=\sqrt{g(x)^{ij}\xi_i\xi_j}$. As a result,
$$
\langle \nu_1(\xi_1),\nu_2(\xi_2) \rangle= \frac{1}{\sqrt{1+4|\xi_1|_{g(x)}^2}\sqrt{1+4|\frac{h}{m}\xi_2|_{g(x)}^2}}+O(\frac{h}{m})
$$

Since $|\xi_1|\gtrsim 1$ and $|\xi_2| \lesssim 1$\footnote{Without loss of generality, we can assume that $||g^{ij}-\delta^{ij}||\leq \frac{1}{C}$ for some large enough $C$ on the coordinate patch considered. This  is enough to have  $|\xi|_{g(x)}\sim |\xi|$.}, we get that $\eqref{transversality schrodinger}$ holds true if $\frac{h}{m}$ is small enough.

The proof of $\eqref{BLS on M T}$ follows by splitting the time interval $[0,T]$ into pieces of length $N_1^{-1}$. That of $\eqref{BLS on M T=1}$ follows by setting $T=1$ in $\eqref{BLS on M T=1}$ when $N_1 \geq 1$ and by using the $L_t^\infty L_x^2$ estimates and S\"obolev's inequality if $N_1\leq 1$.
\endproof

\remark If $P(D)$ is a differential operator on $M$ of degree $n$, then $P(D)e^{iht\Delta} u_0$ has the following expression:

$$
P(D)e^{iht\Delta}u_0(x)=h^{-n} \tilde T'_{h}u_0(t,x)+R'_{h}u_0(t,x)
$$

where $\tilde T '_{h}$ and $R'_{h}$ are operators of the same form as $T_{h}$ and $R_{h}$. In particular, $T'_{h}$ has an expression as in $\eqref{def of tilde T_h}$ (just with different $a$) and $R'_{h}$ obeys similar estimates to $\eqref{Remainder terms}$ (by choosing $h$ small enough). Similar expressions for $e^{imt\Delta}v_0$ allow us, using the exact same analysis performed above, to get:

\begin{corollary}\label{Strichartz on rescaled manifold with differential operators}
Suppose the $u_0, v_0 \in L^2(M)$ are spectrally localized around $N_1,N_2 \in 2^{\Z}$ respectively as in Corollary \ref{Strichartz on rescaled manifold}. Let $P(D)$ and $Q(D)$ be differential operators on $M$ of orders $n$ and $m$ respectively:

\begin{equation}\label{bilinear estimate with differential operators}
||P(D)e^{it\Delta} u_0 Q(D)e^{it\Delta}v_0||_{L^2([0,T]\times M)} \leq N_1^n N_2^m \Lambda(T,N_1,N_2)||u_0||_{L^2(M)}||v_0||_{L^2(M)}
\end{equation}

where $\Lambda(T,N_1,N_2)$ is given in $\eqref{Lambda T}$.
\end{corollary}

This variant will be useful in some applications of the bilinear Strichartz estimates proved here (see \cite{H} for example).

\section{Further Results and Remarks}\label{further results}

\subsection{Bilinear Inhomogeneous Estimates:}
Here we will present some inhomogeneous versions of the bilinear estimates proved in the previous section. We will assume that $u(t)$ and $v(t)$ solve the inhomogeneous Schr\"odinger equation with forcing terms $F$ and $G$ respectively. More precisely:

\begin{eqnarray}
i\partial_t u +\Delta u=F\label{forced SE u}\\
i\partial_t v+\Delta v=G \label{forced SE v}.
\end{eqnarray}

$F$ and  $G$ can be assumed to be a priori in $C^\infty$\footnote{This assumption can be removed a posteriori using standard density arguments.}. The question now is to determine estimates for $||u v||_{L_{t,x}^2}$ in terms of the initial data $u(0)=u_0, v(0)=v_0$ and the forcing terms $F$ and $G$. 

We will prove two types of inhomogeneous estimates: one corresponding to spectrally localized functions generalizing $\eqref{BLS on M}$ and another is a time $T=1$ estimate generalizing $\eqref{BLS on M T=1}$. 

\begin{theorem}\label{inhomogeneous BLS}
Suppose $u(t)$ and $v(t)$ solve the inhomogeneous Schr\"odinger equations $\eqref{forced SE u}$ and $\eqref{forced SE v}$ with initial data $u(0)=u_0$ and $v(0)=v_0$ respectively. Also suppose that $(q,r)$ and $(\tilde q, \tilde r)$ are two Schr\"odinger admissible exponents .

\begin{enumerate}
\item If $u(t)=\varphi(\frac{\sqrt{-\Delta}}{N_1})u(t)$ and $v(t)=\varphi(\frac{\sqrt{-\Delta}}{N_2})v(t)$ for all $t$, then

\begin{equation}\label{inhomogeneous BLS dyadic}
||u v||_{L^2_{t,x}([0,\frac{1}{N_1}]\times M)}\lesssim \left(\frac{N_2^{(d-1)/2}}{N_1^{1/2}}\right)\left(||u_0||_{L^2(M)}+||F||_{L_t^{q'}L_x^{r'}}\right)\left(||v_0||_{L^2(M)}+||G||_{L_t^{\tilde q'}L_x^{\tilde r'}}\right),
\end{equation}

where for any $p\in [1,\infty]$, $p'$ denotes its conjugate exponent $\frac{1}{p}+ \frac{1}{p'}=1$.

\item In general, for any $\delta>0$ we have:

\begin{equation}\label{inhomogeneous BLS T=1}
||u v||_{L^2_{t,x}([0,1]\times M)}\lesssim \left(||u_0||_{H^\delta(M)}+||(\sqrt{1-\Delta})^{\delta+\frac{1}{q}} F||_{L_t^{q'}L_x^{r'}}\right)\left(||v_0||_{H^{1/2-\delta}(M)}+||(\sqrt{1-\Delta})^{1/2-\delta+\frac{1}{\tilde q}}G||_{L_t^{\tilde q'}L_x^{\tilde r'}}\right).
\end{equation}

\end{enumerate}

\end{theorem}

For the proof, we will need the Christ-Kiselev lemma \cite{CK} which we state following Smith and Sogge in \cite{SS}:

\begin{lemma}\label{Christ Kiselev}
Let $X$ and $Y$ be Banach spaces and $K(t,x)$ a continuous function taking values in $B(X,Y)$, the space of bounded linear mappings from $X$ to $Y$. Suppose that $-\infty \leq a <b \leq \infty$ and let

$$
Tf(t)=\int_a^b K(t,s) f(s)ds.
$$

Suppose that
$$
||Tf||_{L^q([a,b];Y)}\leq C||f||_{L^p([a,b];X)},
$$

and define the lower triangular operator 

$$
Wf(t)=\int_a^t K(t,s) f(s) ds.
$$

Then if $1\leq p < q\leq \infty$:

$$
||Wf||_{L^q([a,b];Y)}\lesssim C||f||_{L^p([a,b];X)}.
$$
\end{lemma}

\emph{Proof of Theorem \ref{inhomogeneous BLS}:} We start by proving the spectrally localized version in $\eqref{inhomogeneous BLS dyadic}$. The integral equations satisfied by $u(t)$
 and $v(t)$ are given by Duhamel's formula:
 
 $$
 u(t)=e^{it\Delta}u_0 -i\int_0^t e^{i(t-s)\Delta} F(s)ds,\;\; \; \; v(t)=e^{it\Delta}v_0 -i\int_0^t e^{i(t-s)\Delta}G(s) ds.
 $$

As a result, 
\begin{equation}\label{product u(t)v(t)}
\begin{split}
u(t)v(t)=e^{it\Delta}u_0 e^{it\Delta}v_0 -ie^{it\Delta}u_0\int_0^t e^{i(t-s)\Delta}G(s)ds -ie^{it\Delta}v_0\int_0^t e^{i(t-s)\Delta}F(s)ds\\
- \int_0^t e^{i(t-s)\Delta}F(s)ds \int_0^t e^{i(t-r)\Delta}G(r)dr.
\end{split}
\end{equation}

Recall that $u_0, u(t), F(t)$ are all spectrally localized at dyadic scale $N_1$ and $v_0, v(t), G(t)$ localized at scale $N_2$. The estimate for the first term on the RHS of $\eqref{product u(t)v(t)}$ is the bilinear Strichartz estimate proved in the previous section. We turn to the second term. Applying the Christ-Kiselev lemma (with $Y=L_t^{\tilde q '}L_x^{\tilde r'}$, $X=L^2_{t,x}([0,\frac{1}{N_1}]\times M)$, and $C\sim \frac{N_2^{(d-1)/2}}{N_1^{1/2}}||u_0||_{L^2(M)}$), it is enough to show that:

$$
||e^{it\Delta} u_0 \int_0^{1/N_1}e^{i(t-s)\Delta}G(s)ds||_{L^2_{t,x}([0,\frac{1}{N_1}]\times M)}\lesssim \frac{N_2^{(d-1)/2}}{N_1^{1/2}}||u_0||_{L^2(M)}||G||_{L_t^{\tilde q '}L_x^{\tilde r'}}.
$$

But this follows from the bilinear estimate $\eqref{BLS on M}$ and 

$$
||\int_0^{1/N_1} e^{-is\Delta}\varphi(\frac{\sqrt{-\Delta}}{ N_1})G(s)ds||_{L_x^2(M)}\lesssim ||G||_{L_t^{\tilde q '}L_x^{\tilde r'}}
$$

which is the dual estimate to $\eqref{semiclassical linear estimate}$.

The third term on the RHS of $\eqref{product u(t)v(t)}$ is estimated similarly. For the fourth term, we first apply the Christ-Kiselev lemma to reduce the estimate to the following:

\begin{align*}
||\int_0^{1/{N_1}} e^{i(t-s)\Delta}F(s)ds \int_0^t e^{i(t-r)\Delta}G(r)dr||_{L^2_{t,x}([0,\frac{1}{N_1}]\times M)}
\\
=||e^{it\Delta}\left(\int_0^{1/{N_1}} e^{-is\Delta}F(s)ds \right)\int_0^t e^{i(t-r)\Delta}G(r)dr||_{L^2_{t,x}}\\
\lesssim \frac{N_2^{(d-1)/2}}{N_1}||\int_0^{N_1^{-1}}e^{-is\Delta}F(s)ds||_{L^2(M)}||G||_{L_t^{\tilde q '}L_x^{\tilde r'}}\\
\lesssim ||F||_{L_t^{q'}L_x^{r'}}||G||_{L_t^{\tilde q '}L_x^{\tilde r'}}
\end{align*}

where in the first inequality we apply the same analysis as that used to estimate the second and third term on the RHS of $\eqref{product u(t)v(t)}$ (or apply Christ-Kiselev lemma again) while in the second we use the dual homogeneous Strichartz estimate. This finishes the proof of $\eqref{inhomogeneous BLS dyadic}$.

We now turn to the time $1$ estimate $\eqref{inhomogeneous BLS T=1}$. We start by mentioning that the first term on the RHS of $\eqref{product u(t)v(t)}$ satisfies the needed estimate:

$$
||e^{it\Delta}u_0 e^{it\Delta} v_0||_{L^2([0,1]\times M)}\lesssim ||u_0||_{H^\delta}||v_0||_{H^{1/2-\delta}}.
$$

This follows directly by splitting into Littlewood-Paley pieces: $u=\sum_{\substack{N_1\geq 1\\(dyadic)}}u_{N_1}$ and $v=\sum_{\substack{N_2\geq 1 \\ (dyadic)}}v_{N_2}$ and estimating as follows:

\begin{align*}
||e^{it\Delta}u_0 e^{it\Delta}v_0||_{L^2_{t,x}([0,1]\times M)}\leq& \sum_{N_1\leq N_2 }||e^{it\Delta}u_{N_1} e^{it\Delta}v_{N_2}||_{L^2_{t,x}} +\sum_{N_1>N_2}||e^{it\Delta}u_{N_1} e^{it\Delta}v_{N_2}||_{L^2_{t,x}}\\
\lesssim& \sum_{N_1 \leq N_2}N_1^{(d-1)/2}||u_{N_1}||_{L^2}||v_{N_2}||_{L^2}+\sum_{N_2<N_1}N_2^{(d-1)/2}||u_{N_1}||_{L^2}||v_{N_2}||_{L^2}\\
\lesssim& \sum_{N_1 \leq N_2} \frac{N_1^{(d-1)/2 -\delta}}{N_2^{(d-1)/2-\delta}}||u_{N_1}||_{H^\delta}||v_{N_2}||_{H^{(d-1)/2-\delta}}\\
&+\sum_{N_2<N_1}\frac{N_2^\delta}{N_1^\delta}||u_{N_1}||_{H^\delta}||u_{N_2}||_{H^{(d-1)/2-\delta}}\\
\lesssim& ||u||_{H^\delta}||v||_{H^{(d-1)/2-\delta}}
\end{align*}

where we have used Schur's test to sum in the last step. The rest of the proof of $\eqref{inhomogeneous BLS T=1}$ follows as that of $\eqref{inhomogeneous BLS dyadic}$ above except that here we use the estimate dual to $\eqref{BGT linear estimate}$ given by:

$$
||\int_0^1e^{i(t-s)\Delta}F(s)ds||_{L^2(M)}\lesssim ||(\sqrt{1-\Delta})^{\frac{1}{q}}F||_{L_t^{q'}L_x^{r'}([0,1]\times M)}.
$$

\endproof

\subsection{Bilinear Estimates of mixed type:} Here we present an instance of a mixed-type bilinear estimate of Schr\"odinger-wave type that can be proved using Theorem \ref{Bilinear FIO estimate}. Constant coefficient versions of such estimates are often useful when studying coupled Schr\"odinger-wave systems such as the Zakharov system (see \cite{BHHT} for instance). Theorem \ref{mixed BLS} below serves as an example of a variable coefficient Schr\"odinger-wave bilinear estimates and has potential applications in studying Zakharov systems (or other Schr\"odinger-wave systems) on manifolds.

\begin{theorem}\label{mixed BLS}
Suppose $u_0, v_0\in L^2(M^d)$ are spectrally localized at dyadic scales $N_1$ and $N_2$ as above with $1 \ll N_1$. Then the following estimate holds:

\begin{equation}\label{Schr-Wave BLS}
\left|\left|e^{it\Delta}u_0 e^{\pm it|\nabla|}v_0\right|\right|_{L^2_{t,x}([-\frac{1}{N_1},\frac{1}{N_1}]\times M)} \lesssim_M \frac{\min(N_1,N_2)^{\frac{d-1}{2}}}{N_1^{1/2}}||u_0||_{L^2(M)}||v||_{L^2(M)}.
\end{equation}

\end{theorem}

Of course, an estimate over the time interval $[0,T]$ follows as well by splitting into into pieces of length $\frac{1}{N_1}$.

\proof We present the proof in the case of the forward half wave operator, the proof for the backwards operator being similar. As before, we use the parametrix for $e^{it|\nabla|} v_0$ which is given, up to a smoothing remainder $R_m v_0$, by the oscillatory integral:

$$
S_m^W v_0=\frac{1}{(2\pi m)^d}\int_{\R^d}e^{\frac{i}{m}\psi(t,x,\xi_2)}a(t,x,\xi_2)\widehat{\tilde v_0}(\frac{\xi_2}{m})d\xi_2
$$

where $\psi$ is a non-degenerate phase function (in particular $\det \left(\frac{\partial^2}{\partial \xi \partial x}\tilde \psi\right)\neq 0$) and homogeneous in $\xi_2$ of degree 1 and $\tilde v_0$ is a microlocalization of $v_0$ as explained in section \ref{proof of BLS on M} (cf. \cite{H2} Chapter XXIX). As before, we used the convention that $h=\frac{1}{N_1}$ and $m=\frac{1}{N_2}$. As a result, we have:

$$
\left|\left|e^{it\Delta}u_0 e^{it|\nabla|}v_0\right|\right|_{L^2_{t,x}([-\frac{\alpha}{N_1},\frac{\alpha}{N_1}]\times M)} =
h^{1/2}\left|\left|e^{iht\Delta}u_0 e^{iht|\nabla|}v_0\right|\right|_{L^2_{t,x}([-\alpha,\alpha]\times M)}. 
$$

Ignoring the smooth remainder terms $R_h$ and $R_m$ (as they are inconsequential as in section \ref{proof of BLS on M}) we get that $\eqref{Schr-Wave BLS}$ follows from the estimate:

\begin{align*}
||\tilde T_h u_0(t,x) \tilde S_m^W v_0(ht,x)||_{L^2_{t,x}([-\alpha,\alpha]\times \R^d)}\lesssim& \frac{1}{(hm)^{d/2}}\min(m,h)^{d/2}\max(m,h)^{1/2}||\tilde u_0||_{L^2(\R^d)}||\tilde v_0||_{L^2(\R^d)}\\
=&C \max(m,h)^{-(d-1)/2}||\tilde u_0||_{L^2(\R^d)}||\tilde v_0||_{L^2(\R^d)}
\end{align*}

This inequality follows by applying Theorem \ref{bilinear estimate} with the non-degenerate phase functions $\phi(t,x,\xi_1)=\tilde \phi(t,x,\xi_1)$ and $\psi(t,x,\xi_2)=\tilde \psi(ht,x,\xi_2)$. The transversality condition is directly verified as follows:
the normal vectors to the two surfaces:

\begin{align*}
S_\phi&: \xi_1\mapsto \nabla_{t,x}\phi(t,x,\xi_1)=(\nabla_x \tilde \phi(t,x,\xi_1), \partial_t \tilde \phi(t,x,\xi_1)) \\
S_{\psi} &:\xi_2 \mapsto \nabla_{t,x}\psi(t,x,\xi_2)=(\nabla_x \tilde \psi(ht,x,\xi_2), h\partial_t \tilde \psi(ht,x,\xi_2))
\end{align*}

can be written as $\nu_1=(\eta_1, \tau_1)$ and $\nu_2=(\eta_2, \tau_2)$ with $\eta_1, \eta_2 \in \R^n$ and $\tau_1, \tau_2 \in \R$. The fact that $\langle \nu_2, \frac{\partial^2}{\partial \xi \partial (x,t)}\psi\rangle =\vec{0}$ implies that $\langle \eta_2, \frac{\partial^2}{\partial \xi \partial x}\tilde \psi(ht,x,\xi_2) \rangle+h\tau_2\partial_t \partial _\xi \tilde \psi(ht,x,\xi_2)=\vec{0}$ which implies that 

$$\eta_2=-h\tau_2\langle \partial_t \partial _\xi \tilde \psi,\left[\frac{\partial^2}{\partial \xi \partial x}\tilde \psi\right]^{-1} \rangle=O(h).$$

This gives that 

$$
\langle \nu_1,\nu_2 \rangle \leq |\tau_1 \tau_2|+O(h)\leq |\tau_1|+O(h).
$$

As a result, the transversality condition $\eqref{transversality condition}$ holds if $h\ll 1$ (i.e.~$N_1 \gg 1$) and $|\tau_1| <1$ which is the case since  $\tau_1=\frac{-1}{\sqrt{1+4|\xi|_{g(x)}^2}}$ and $|\xi_1|\gtrsim 1$ (see end of the proof of Theorem \ref{BLS on M theorem}).

\endproof

\subsection{Applications in PDE} The bilinear estimate $\eqref{BLS on M T=1}$ directly implies local well-posedness for 2-dimensional cubic NLS:

\begin{equation}\label{cubic NLS}
\begin{split}
i\partial_t u +\Delta u&=|u|^2u\\
u(t=0)&=u_0\in H^s(M^2)
\end{split}
\end{equation}

in $X^{s,b} \subset C_tH_x^s$ spaces for all $s>1/2$ and some $b>1/2$. It should be noted that local well-posedness of $\eqref{cubic NLS}$ in $C_tH^s$ for $s>1/2$ has already been proven in \cite{BGT} using linear Strichartz estimates. Here $X^{s,b}$ is the closure of $C_0^\infty(\R \times M)$ in the norm:

$$
||u||_{X^{s,b}}=\left(\int_\R \sum_{\nu} \langle \tau+\nu \rangle^{2b}\langle \nu \rangle^s ||\widehat{\pi_{\nu}u}(\tau)||_{L^2(M)}^2d\tau \right)^{1/2}
$$

where the sum runs over the distinct eigenvalues of the Laplacian and $\pi_\nu$ is the projection onto the eigenspace corresponding to the eigenvalue $\nu$. It is worth remarking that $\eqref{BLS on M}$ translates into the following estimate for functions $u, v \in C_0^\infty(\R \times M)$ satisfying $u(t)=\mathbf{1}_{[N_1,2N_1)}(\sqrt{-\Delta}) u(t)$  and $v(t)=\mathbf{1}_{[N_2,2N_2)}(\sqrt{-\Delta}) v(t)$:

\begin{equation}\label{bilinear Xsb}
||u v||_{L^2(\R\times M)} \lesssim \min(N_2,N_1)^{1/2}||u||_{X^{0,b}}||v||_{X^{0,b}}
\end{equation}

for any $b>1/2$ (cf \cite{BEE},\cite{H}). Using this and a standard dyadic decomposition one can prove the crucial cubic estimate that yields local well-posedness via Picard iteration (see \cite{BEE} for example).

One interesting application of Theorem \ref{BLS on M theorem} is that of proving global well-posedness of $\eqref{cubic NLS}$ for $s<1$. As mentioned in the introduction, the bilinear Strichartz estimate $\eqref{BLS on M T}$ on the time interval $[0,T]$ translates into a bilinear Strichartz estimate on the rescaled manifold $\lambda M$ over the time interval $[0,1]$. Here $\lambda M$ can either be viewed as the Riemmannian manifold $(M, \frac{1}{\lambda^2}g)$ or by embedding $M$ into some ambient space $R^N$ and then applying a dilation by $\lambda$ to get $\lambda M$. The relevant result was cited in the introduction in corollary \ref{Strichartz on rescaled manifold}: if  $u_0,v_0 \in L^2(\lambda M)$ are spectrally localized around $N_1$ and $N_2$ respectively, with $N_2\leq N_1$. Then

\begin{align*}
||e^{it\Delta_{\lambda}} u_0 e^{it\Delta_{\lambda}}v_0||_{L^2([0,1]\times \lambda M)}
\lesssim& \Lambda(\lambda^{-2},\lambda N_1,\lambda N_2)||u_0||_{L^2(\lambda M)}||v_0||_{L^2(\lambda M)}\\
\lesssim& \left\{
    \begin{array}{ll}
        \left(\frac{N_2}{N_1}\right)^{1/2} ||u_0||_{L^2(\lambda M)}||v_0||_{L^2(\lambda M)} & \mbox{if } \lambda \gg N_1 \\
       \left(\frac{N_2}{\lambda}\right)^{1/2}||u_0||_{L^2(\lambda M)}||v_0||_{L^2(\lambda M)} & \mbox{if } \lambda \lesssim N_1.
    \end{array}
\right.
\end{align*}

This estimate turns out to be crucial in \cite{H} where it is proved that $\eqref{cubic NLS}$ is globally well-posed for all $s>2/3$. This generalizes, \emph{without any loss in regularity}, a similar result of Bourgain \cite{B4}(see also \cite{dSPST}) where global well-posedness for $s>2/3$ is proved for the torus $\T^2$. Global well-posedness for $s\geq 1$ follows using conservation of energy and standard arguments. To go below the energy regularity $s=1$, the I-method of Colliander, Keel, Staffilani, Takaoka, and Tao should be used and most of the analysis is done on $\lambda M$ rather than $M$. As a result, the factor of $\frac{1}{\lambda^{1/2}}$ on the R.H.S. of $\eqref{the N_2/lambda decay}$ in the range $\lambda \lesssim N_1$ becomes crucial to get the full regularity range of $s>2/3$ (see \cite{H}).

\end{document}